\documentclass[11pt]{article}%
\usepackage{graphicx}
\usepackage{subcaption}
\usepackage{makeidx}
\usepackage{amssymb}
\usepackage[all]{xy}
\usepackage{float}
\usepackage[latin1]{inputenc}
\usepackage{amsfonts}
\usepackage{amsmath}
\usepackage{amssymb}
\usepackage{url}
\usepackage{hyperref}
\usepackage{graphicx}%
\usepackage{mathabx}
\usepackage{wrapfig}

\usepackage{amsthm}
\usepackage{url}
\usepackage{hyperref}
\providecommand{\U}[1]{\protect\rule{.1in}{.1in}}

\usepackage{xcolor}

\newtheorem{prop}{Proposition}[section]
\newtheorem{cor}[prop]{Corollary}
\newtheorem{defi}[prop]{Definition}

\newtheorem{lem}[prop]{Lemma}

\newtheorem{theo}[prop]{Theorem}
\newtheorem{rem}[prop]{Remark}

\def\d{\delta}

\newcommand{\tr}{\mbox{\rm Tr}}

\newcommand{\CC}{\mathbb{C}}

\newcommand{\EE}{\mathbb{E}}
\newcommand{\FF}{\mathbb{F}}
\newcommand{\HH}{\mathbb{H}}

\newcommand{\LL}{\mathbb{L}}

\newcommand{\NN}{\mathbb{N}}
\newcommand{\PP}{\mathbb{P}}

\newcommand{\RR}{\mathbb{R}}

\newcommand{\Ca}{ {\cal C }}

\newcommand{\La}{ {\cal L }}
\newcommand{\Na}{ {\cal N }}
\newcommand{\Ka}{ {\cal K }}

\newcommand{\Ea}{ {\cal E }}
\newcommand{\Sa}{ {\cal S }}

\newcommand{\Ga}{ {\cal G }}

\newcommand{\Ma}{ {\cal M }}

\newcommand{\Ha}{ {\cal H }}
\newcommand{\Ja}{ {\cal J }}
\newcommand{\Pa}{ {\cal P }}
\newcommand{\Za}{ {\cal Z }}

\newcommand{\Wa}{ {\cal W }}

\newcommand{\point}{\mbox{\LARGE .}}

\newcommand{\cqfd}{\hfill\blbx \\}
\def\blbx{\hbox{\vrule height 5pt width 5pt depth 0pt}\medskip}

\def \PP{\mathbb{P}}
\def \RR{\mathbb{R}}
\def \SS{\mathbb{S}}
\def \EE{\mathbb{E}}

\def \CC{\mathbb{C}}
\def \LL{\mathbb{L}}

\def \WW{\mathbb{W}}
\def \BB{\mathbb{B}}

\begin{document}

  \title{Quantum harmonic oscillators and Feynman-Kac path integrals for
   linear diffusive particles}
  \author{P. Del Moral, E. Horton}


\maketitle

\begin{abstract}
We propose a new solvable class of  multidimensional quantum harmonic oscillators for a
  linear diffusive particle and a quadratic energy absorbing well associated with a semi-definite positive matrix force. Under natural and easily checked controllability conditions,  the ground state  and the zero-point energy are explicitly computed in terms of a positive fixed point of a continuous time algebraic Riccati matrix equation. We also present an explicit solution of normalized and time dependent Feynman-Kac measures in terms of a time varying linear dynamical system coupled with a differential Riccati matrix equation.

   A refined non asymptotic analysis of the stability of these models is developed based on a recently developed Floquet-type representation of time varying exponential semigroups of Riccati matrices. We provide explicit and non asymptotic 
estimates of the exponential decays to equilibrium of Feynman-Kac semigroups in terms of Wasserstein distances or  Boltzmann-relative entropy. 

For reversible models we develop
a series of functional inequalities including de Bruijn identity,
Fisher's information decays, log-Sobolev inequalities, and entropy contraction estimates.
In this context, we also provide a complete and explicit description of all the spectrum and the excited states of the Hamiltonian,  yielding what seems to be the first result of this type for this class of models.

  We illustrate these formulae with the traditional harmonic oscillator associated with
 real time  Brownian particles and Mehler's formula. The analysis developed in this article can also be extended to solve time dependent Schrodinger equations equipped with time varying linear diffusions and quadratic potential functions.\\

\newline


\end{abstract}


\section{Introduction}

\subsection{Description of the models}
 Given $r\geq 1$ and some square $(r\times r)$-matrices $A,R,S$ with real entries,  let 
$\Ha$ be the Hamiltonian differential operator given by the formula
\begin{equation}\label{def-H-intro}
\Ha=-\La+V\quad \mbox{\rm with the potential energy}\quad V(x):=\frac{1}{2}\sum_{1\leq k,l\leq r} x_k~S_{k,l}~x_l
\end{equation}
In the above display, $\La$ stands for the second order differential kinetic energy operator
\begin{equation}\label{def-La-intro}
\La:=\sum_{1\leq k,l\leq r} A_{k,l}~x_l~\partial_{x_k}+\frac{1}{2}\sum_{1\leq k,l\leq r} R_{k,l}~\partial_{x_k,x_l}
\end{equation}
In the present article we shall assume  that $R$ and $S$ are positive semi-definite matrices
and the pairs of matrices $(A,R^{1/2})$ and  $(A^{\prime},S^{1/2})$ are both controllable, 
 in the sense that the $(r\times r^2)$-matrices
\begin{equation}\label{def-contr-obs}
\left[R^{1/2},AR^{1/2}\ldots, A^{r-1}R^{1/2}\right]\quad
\mbox{\rm and}\quad
\left[S^{1/2},A^{\prime}S^{1/2}\ldots, (A^{\prime})^{r-1}S^{1/2}\right] \quad\mbox{have rank $r$}.
\end{equation}
In the above display, $A^{\prime}$ stands for the transposition of the matrix $A$.
Note that the above condition holds trivially when $R$ and $S$ are positive matrices.

The time dependent Schr\"odinger equation and the imaginary time version associated with the hamiltonian $\Ha$ are given respectively by the equations
\begin{equation}\label{def-Schrod-intro}
i\,\partial_t \Psi_t(x)=\Ha(\Psi_t)(x)\quad \mbox{\rm and}\quad
-\partial _t\psi_t(x)=\Ha(\psi_t)
\end{equation}
with prescribed initial conditions.
In the above display, $i\in\CC$ stands for the complex number such that $i^2=-1$. The right hand side equation is obtained by a formal change of time by setting
$\psi_t(x)=\Psi_{-it}(x)$. 
The corresponding evolution equation takes the following form 
\begin{equation}\label{def-Schrod-intro-imaginary}
\partial_t \psi_t(x)=\La(\psi_t)(x)-V(x)\psi_t(x)
\end{equation}
The first term $\La$ represents the generator of a free linear diffusion process $X_t$ with drift matrix and diffusion matrix $R$. The stochastic differential equation associated with this free evolution process is described in more details  in (\ref{lin-Gaussian-diffusion-filtering}). We emphasize that $A$ is not required to be a stable Hurwitz matrix so that $X_t$ can be a transient diffusion process that evolves exponentially fast to $\infty$. As we shall see in the further development of the article,
the controllability condition (\ref{def-contr-obs}) ensures that the trapping force of the potential energy always compensates the delocalization kinetic energy of the diffusion.

For a  twice differentiable function $\psi_0$, the solution of (\ref{def-Schrod-intro-imaginary}) is given by the Feynman-Kac path integral formula
\begin{eqnarray}
\psi_t(x)&=&\Ka_t(\psi_0)(x):=\int\Ka_t(x,dy)~\psi_0(y)\label{def-Ka}\\
&=&\EE\left(\psi_0(X_t)~\exp{\left(-\int_0^t V(X_u)~ds\right)}~|~X_0=x\right)\nonumber.
\end{eqnarray}

The integral operator $\Ka_t$ is sometimes called the Feynman-Kac propagator. Besides its mathematical elegance, the  conditional expectations in (\ref{def-Ka}) can rarely be solved analytically. Moreover, numerical solutions for general diffusions $X_t$ and/or potential functions $V$ that are not necessarily quadratic require extensive calculations, see for instance~\cite{caffarel,cances-tony,dm-04,dm-13} and the references therein. 

There exists a rich literature on the micro-local analysis~\cite{grigis,ivrii,sjostrand} and the semi-classical analysis~\cite{dimasi,Martinez,zworski} of self-adjoint Hamiltonian operators for general smooth potentials and Brownian particle free motions. These powerful mathematical tools provide a precise spectral asymptotic analysis by connecting the Schr\"odinger equation with the classical mechanics of point particles (a.k.a. Bohr correspondence principle) when the diffusion Planck constant tends to $0$. 

Non-asymptotic estimates for general models are rarely studied in the literature and they often rely on proving the existence of limiting unknown mathematical objects  such as quasi-invariant measures, the zero-point energy and the corresponding ground state, see for instance~\cite{dm-04,dm-2000,dm-sch,dm-sch-2} and the more recent articles of Champagnat and Villemonais~\cite{champagnat-1,champagnat-2,champagnat-3,champagnat-4,dm-sch}.

By an elementary second order Taylor expansion, any smooth potential $V$ can be approximated by a harmonic quadratic potential of the form (\ref{def-H-intro}) at the vicinity of a stable equilibrium point. 
The quantum harmonic oscillator corresponding to the  case $A=0$ and diagonal matrices $(R,S)$  is one of the most important quantum-mechanical Hamiltonian systems for which an exact analytical solution is known. To the best of our knowledge, the case $A\not=0$ and non diagonal matrices $(R,S)$ has not been considered in the literature on this subject.
This article provides an analytical solution and a complete theoretical analysis the multivariate quantum harmonic oscillator and related particle absorption processes for general hamiltonian operators given by (\ref{def-H-intro}).\\

A brief description of the main objectives and the main results of this article is provided below:

\begin{itemize}
\item One of the main objective of this article is to construct an explicit closed form solution of the time dependent Schr\"odinger equation (\ref{def-Schrod-intro}) for abstract and general Hamiltonian operators of the form (\ref{def-H-intro}) in the reversible situation, that is when  $R>0$ and $AR=RA^{\prime}$ (see Theorem~\ref{theo-Ka-spec}). In this context we provide a complete description of the entire spectrum of the Hamiltonian operator, including all the excited states in terms of the matrices $(A,R,S)$.
In the non-reversible case, we also provide an explicit description of the zero-point energy  and the ground state of the Hamiltonian in terms of the positive fixed point of a continuous algebraic Riccati matrix equation (a.k.a. CARE, see for instance (\ref{def-h}) and Theorem~\ref{theo-1-intro}). 

\item When the process is not necessarily reversible, our second main objective is to explicitly compute the time varying distribution flow of survival probabilities (\ref{survival-proba}) including  the distributions of a non-absorbed particle (\ref{def-d-non-abs}); see for instance the Gaussian preserving property (\ref{cv-eta-infty}) and the coupled equations (\ref{ricc-intro}). The distribution of  a non-absorbed diffusion and a particle evolving in the ground state (a.k.a. $h$-process) are connected to each other by a Boltzmann-Gibbs transformation (a.k.a. Bayes' rule or Doob's $h$-transform, 
see for instance~(\ref{BG-h0-eta-intro}) and (\ref{Ka-Kah}), as well as Theorem~\ref{theo-c-h0} in the context of path space models). For any initial Gaussian state we show that a non-absorbed particle remains distributed according to a Gaussian probability with a mean vector satisfying a coupled time-varying linear system depending of the solution of a time dependent Riccati differential equation (\ref{ricc-intro}).

\item An important part of the article is concerned with the long time behavior of $h$-processes (\ref{def-X-h-d-intro}) and normalised Feynman-Kac measures (cf. section~\ref{sec-norm-measures} and section~\ref{sec-norm-propagators}), including the convergence of the conditional distribution of a non absorbed particle (\ref{def-d-non-abs}) towards the unique fixed point (a.k.a. quasi-invariant distribution) of a nonlinear semigroup in distribution space (\ref{cv-h0-proba-fixed point}). In the reversible case, the density of these limiting distributions with respect to the reversible measure of the free particle coincides with the ground state of the Hamiltonian (\ref{cv-h0-proba}). 

Our main contributions to the stability analysis of $h$-processes and normalised Feynman-Kac semigroups are twofold:
\begin{itemize}
\item Firstly, we provide explicit and non-asymptotic 
estimates of the exponential decays to equilibrium in terms of Wasserstein distances or  Boltzmann-relative entropy. These results  are summarised in Theorems~\ref{theo-lamda0-K} and~\ref{theo-intro-entrop-wasserstein}. We emphasise that these theorems are valid for non-necessarily reversible models, even when the drift matrix $A$ is unstable, yielding what seems to be the first result of this type for this class of Feynman-Kac particle absorption models.

\item In the reversible case, we analyze the stability properties of the $h$-process with
a series of functional inequalities including de Bruijn identity (\ref{de-Bruijn}),
Fisher's information decays, log-Sobolev inequalities, and entropy contraction estimates (cf. Theorem~\ref{theo-entropy-fisher}). We also deduce Poincar\'e inequality and variance-type exponential decays directly from the spectral theorem~\ref{spectral-theo-intro} (see also Corollary~\ref{Poincare-cor}).

\end{itemize}
\item Last but not least, section~\ref{mcKean-sec} discusses several classes of McKean-Vlasov interpretations of the distribution of a non-absorbed particle.  These probabilistic models and their mean field simulation are defined in terms of a nonlinear Markov process that depends on the distribution of the random states so that the flow of distributions of all random states coincides with the conditional distribution of a non-absorbed particle. For a more thorough discussion on these nonlinear sampling methodologies we refer to the books~\cite{dm-04,dm-13,dm-2000,dmpenev} and the references therein.

 Section~\ref{IPS-section} is dedicated to interacting jump interpretations. Their mean field interpretations coincides with conventional Quantum Monte Carlo methods currently used in numerical physics. Their path space version can be interpreted as the genealogical tree
 associated with the killing and the birth/duplication/selection of walkers. An alternative way of sampling the trajectories of a non absorbed particle backward is provided in
 section~\ref{section-backward} (see for instance Theorem~\ref{theo-backward}).
 
 In section~\ref{EnKF-sec} we present a new class of mean field samplers based on Ensemble Kalman filters and the novel feedback particle filter methodology introduced by Mehta and Meyn and their co-authors in a series of seminal articles~\cite{prashant-1,prashant-2,prashant-3,prashant-4,prashant-5}. To the best of our knowledge, this class of advanced Monte Carlo methods have not been used to solved ground state energies nor to sample non-absorbed particle processes.

\end{itemize}

A more formal discussion on the probabilistic models and the main results presented in this article 
is provided in section~\ref{ss1} and section~\ref{intro-absorption-sec}.  The detailed statements of the main theorems  are presented  in section~\ref{sec-main-results}.

To facilitate the interpretation of the theoretical and numerical physics in the measure theoretical framework used in this article, we end this introduction with some comments on the probabilistic setting. In theoretical and mathematical physics, the Feynman-Kac propagator defined by the integral operator (\ref{def-Ka}) is sometimes written in terms of the exponential of the Hamiltonian  operator with the exponential-type symbol
 $$
\Ka_t:=e^{-t\Ha}\quad \mbox{\rm or in the bra-kets formalism}\quad \Ka_t(\psi_0)=\vert e^{-t\Ha}\vert\psi_0\rangle.
$$ 
The exponential notation is compatible with finite space models and the matrix notation of 
the continuous one-parameter semigroup for time homogenous models. The bra-ket notation (a.k.a. Dirac notation) is also used to represents linear projection forms acting on Hilbert spaces associated with some reversible or some stationary measure, such as the Lebesque measure for the harmonic oscillator.
 
The present article deals with different types of non necessarily stationary stochastic processes,
including the  free evolution process $X_t$ discussed in (\ref{def-Ka}), $h$-processes and McKean-Vlasov jump or diffusion-type processes.  
Apart in the reversible situation in which spectral theorems are stated on the Hilbert space associated with a reversible measure, the use of the exponential symbol or the use of 
the bra-kets formalism  is clearly not adapted to  represent different  expectations with respect to different types of stochastic and non-necessarily reversible processes. 

To analyze these general stochastic models, we have chosen to only use elementary and standard measure theory notation such as (\ref{def-Ka}). The integral actions of a given integral operator $K_t(x,dy)$ on the right on functions $f(y)$ (such as (\ref{def-Ka})) and on the left  on measures $\mu(dx)$ (such as (\ref{cv-h0-proba-fixed point})) are clearly compatible with finite space models and matrix notation. The left action $\mu\mapsto\mu K_t$ maps measures into measures, while the right action $f\mapsto K_t(f)$ maps functions into functions 
$$
(\mu K_t)(dy):=\int \mu(dx)~K_t(x,dy)\quad \mbox{\rm and}\quad K_t(f)(x):=\int K_t(x,dy)~f(y).
$$
For finite or countable state space models the integrals are clearly replaced by finite or countable sums, the integral operator $K_t(x,dy)$ is replaced by a square matrix, the function $f(y)$ by a column vector, and by duality, the measure $\mu(dx)$ is represented by a row vector. These are the only notation from measure theory used in the present article.

For any $s,t\geq 0$ the integral operators $\Ka_t$ introduced in (\ref{def-Ka}) satisfy the semigroup property
$$
\Ka_{s+t}(x,dz)=(\Ka_s\Ka_t)(x,dy):=\int~\Ka_s(x,dz)~\Ka_t(z,dy)\Longrightarrow\psi_{s+t}=\Ka_{s}(\psi_t).
$$
In terms of left action bra-kets, defining $\mu_{\varphi}(dx):=\varphi(x)dx$, Fubini's theorem yields
$$
\begin{array}[t]{rcl}
\displaystyle \langle \varphi\vert e^{-s\Ha}\vert \psi_t\rangle &=&\displaystyle \int dx~ \varphi(x) ~\Ka_s(x,dy)~\psi_t(dy)= (\mu_{\varphi}\Ka_s)(\psi_t)\\
&&\\
&=&\displaystyle  \mu_{\varphi}((\Ka_s\Ka_t)(\psi_0))=\mu_{\varphi}\Ka_{s+t}(\psi_0)=\langle \varphi\vert e^{-(s+t)\Ha}\vert \psi_0\rangle.
\end{array}
$$

 \subsection{Harmonic oscillator for linear diffusions}\label{ss1}
One of the main questions of quantum mechanics is to find the quantum numbers $n$, the eigenstates $h_n$ and the energy levels  $\lambda_n$  of the Hamiltonian operator introduced in (\ref{def-H-intro}); that is, to find a sequence of functions $h_n(x)$ in some Hilbert space and some energy levels $\lambda_n\in \RR_+:=[0,\infty[$ satisfying for any quantum numbers $n$ the time independent Schr\"odinger equation 
\begin{equation}\label{Psi-spectral-intro-1}
\Ha(h_n)(x)=\lambda_n~h_n(x)~~\Longleftrightarrow~~ \Ka_t(h_n)(x)=\exp{\left(-\lambda_n t\right)}~h_n(x).
\end{equation}

The complete answer to this question is rather well known for the conventional isotropic harmonic oscillator associated with a null matrix $A=0$ and diagonal matrices $(R,S)$.  The one dimensional case with $A=0$ corresponds to the well know harmonic oscillator treated in any textbook in quantum mechanics, see for instance~\cite{davidov,flugge,landau-2,schiff,speck}. In the multidimensional case, the Hamiltonian resumes to the sum of independent operators in each dimension. The resulting energy levels coincide with the tensor product of energy levels in each dimension. The isotropic harmonic oscillator corresponds to the case where $S=\rho I$, for some constant $\rho>0$. The case $A=0$ and non diagonal matrices $S$ arise in the analysis of coupled harmonic oscillators, see for instance~\cite{cook,davidova,Makarov,mcdermott,Park} and references therein. Coupled harmonic oscillators arise in a variety of applications including quantum and nonlinear physics~\cite{fano,paz},
 quantum cryptography and communication~\cite{eker,benett}, quantum teleportation~\cite{samuel}, as well as in biophysics~\cite{romero,halpin} and in molecular chemistry~\cite{ikeda,melan}. 

To the best of our knowledge the case $A\not=0$ has not been considered in quantum mechanics literature, the hypothesis of universal Brownian particle velocities in real time is always in force in all the studies published in this field.
The main objective of this article is to extend conventional 
quantum harmonic oscillators to linear drift-type particle diffusions and general potential functions associated with some quadratic form. This class of models differs from the damped quantum harmonic oscillators with Ornstein-Uhlenbeck stable diffusions  in imaginary time discussed in the series or articles~\cite{aguilar,bateman-2,bouchaud,cordero,cordero-2,dekker,landau}.

For matrices $(A,R,S)$ satisfying the controllability condition (\ref{def-contr-obs}), we provide an explicit description of the zero-point energy $\lambda_0$ and the ground state $h_0$ of the Hamiltonian (see Theorem~\ref{theo-1-intro}) in terms of the positive fixed point  of a continuous time algebraic Riccati matrix equation. More precisely,  we have 
 \begin{equation}\label{def-h}
 \begin{array}{l}
\displaystyle Q_{\infty}>0 \quad \mbox{\rm and}\quad A^{\prime}Q_{\infty}+Q_{\infty}A-Q_{\infty}RQ_{\infty}+S=0\\
\\
\displaystyle\Longrightarrow\quad \lambda_0=\frac{1}{2}\,\tr\left(R\,Q_{\infty}\right)\quad \mbox{\rm and}\quad
h_0(x)=\exp{\left(-\frac{1}{2}~x^{\prime}Q_{\infty}\,x\right)}.
\end{array}
\end{equation}
Here, and in the rest of the article, $\tr(\point)$ stands for the trace operator.
 As a rule in the present article, the state vectors  $x\in \RR^d$ are column vectors and  $x^{\prime}$ stands for the transposed row vector. 
 
 Riccati equations such as the one discussed above play a central role in signal process and optimal control theory, starting with the pioneering work of Kalman in the beginning of the 1960s, see for instance~\cite{abou-kandil03,anderson-moore,Bittanti91,Lancaster1995,reid}, and the more recent articles~\cite{BishopDelMoralMatricRicc,bp-21} in the context of filtering theory. To the best of our knowledge, their application in the context of quantum harmonic oscillators seems to be new. 
 
 For large scale problems, the numerical solving of the algebraic Riccati matrix equation (\ref{def-h}) using exact or inexact Kleinman-Newton type methods  is  generally impractical~\cite{feitzinger,Lancaster1995,Kleinman}. Several improvements have been suggested to avoid the degeneracy of the residuals arising in the Lyapunov recursions associated with these sequential gradient type estimates~\cite{benner-1,benner-2,guo,kelley}. An alternative approach is to use Diffusion Monte Carlo methods and extended versions of Ensemble Kalman type methodologies used in signal processing and information theory.
 
In the reversible situation, that is when $AR=RA^{\prime}$, we solve the Schr\"odinger eigenvalue problem with the imaginary time technique. In this context the entire spectrum of $\Ha$ can be computed explicitly in terms of the matrices $(A,R,S)$. For instance, the ground state can be computed with the formula
\begin{equation}\label{ref-Weyls-2-intro-0}
Q_{\infty}=R^{-1}A+R^{-1}(A^2+RS)^{1/2}.
\end{equation}
In the above display $(A^2+RS)^{1/2}$ stands for the square root that has positive eigenvalues. A proof of the above assertion is provided in section~\ref{rev-section}.

In section~\ref{spectral-intro-sec} we shall see that the energy levels $\lambda_n$ are indexed by multiple index quantum numbers $n=(n_1,\ldots,n_r)\in \NN^r$ and given by the formulae
\begin{equation}\label{first-example}
\lambda_n=\frac{1}{2}~\tr\left(A\right)+\sum_{1\leq i\leq r}\left(n_i+\frac{1}{2}~\right)\sqrt{\vert\lambda_i(A^2+RS)\vert},
\end{equation}
where $\lambda_i(A^2+RS)$ stands for the non-negative eigenvalues of $(A^2+RS)$.
The excited eigenstates $h_n$ with the energy level $\lambda_n$ are defined on the Hilbert space $\LL_2(\mu)$ associated with a locally finite Gibbs-type measure $\mu$ that only depends on the matrices (A,R) (see Theorem~\ref{theo-Ka-spec}). We already mention that $\mu$ is Gaussian if and only if the drift matrix $A$ is Hurwitz.

After decomposing the initial state $\psi_0=\Psi_0$ into the $h_n$-basis discussed above, we apply the time evolution at each energy level $\lambda_n$. Reassembling all eigenstates we obtain the solution of both the time dependent Schr\"odinger equation and the imaginary time version discussed in (\ref{def-Schrod-intro}); that is, we have that
\begin{equation}\label{Psi-spectral-intro}
\Psi_{t}(x)=\sum_{n\in\NN^r} e^{-i\lambda_nt}~h_n(x)~\int~\Psi_0(y)~h_n(y)~\mu(dy)\quad \mbox{\rm and}\quad
\psi_t(x)=\Psi_{-it}(x).
\end{equation}
The null quantum number $(0,\ldots,0)\in \NN^r$ will correspond to the bottom of the spectrum of the Hamiltonian $\Ha$ and to simplify notation we shall write $\lambda_0$ and $h_0$ instead of $\lambda_{(0,\ldots,0)}$ and $h_{(0,\ldots,0)}$. A detailed description of 
the measure $\mu$, the energy levels and the corresponding eigenstates discussed in 
(\ref{Psi-spectral-intro-1}) and (\ref{Psi-spectral-intro}) is provided in the end of section~\ref{rev-models-intro}.

We illustrate these spectral decompositions in section~\ref{sec-illustrations}.
One dimensional models are discussed in section~\ref{sec-one-d}. 
In this situation, it is clear from (\ref{first-example}) that the
the trapping force of the potential energy always compensates the delocalisation kinetic energy of the diffusion even when $A>0$ is very large.

In section~\ref{Mehler-sec} we show how to recover directly Mehler's formula
 from our spectral decompositions. The multidimensional quantum harmonic oscillator discussed in section~\ref{multi-dimn-osc-sec} corresponds to the null drift $A=0$ and diagonal matrices $(R,S)$. 
 
\subsection{Particle absorption processes}\label{intro-absorption-sec}
Consider a  process
$X^c_t$ starting from $X^c_0=X_0$, evolving as the diffusion $X_t$ and killed  with rate $V(X^c_t)$. We denote by $\tau^c$ the random killing time of the process.
In this interpretation, the Feynman-Kac propagator discussed in (\ref{def-Ka}) takes the following form
\begin{equation}\label{FKP-2}
\psi_t(x)=\Ka_t(\psi_0)(x)=\EE(\psi_0(X^c_t)~1_{\tau^c\geq t}~|~X^c_0=x).
\end{equation}
An important question arising in applied probability and rare event analysis is to study the long time behavior of the conditional probability of the process $X^c_t$ with respect to the non-absorption event and starting from a random variable $X^c_0=X_0$ with distribution $\eta_0$.  In a more synthetic form this distribution is given by the formula
\begin{equation}\label{def-d-non-abs}
\eta_t(dx):=\PP_{\eta_0}(X^c_t\in dx~|~\tau^c>t).
\end{equation}
Equivalently, for any bounded measurable function $f$ on $\RR^r$ we have the integral formula
$$
\eta_t(f):=\int f(x)~\eta_t(dx)=\EE(f(X^c_t)~|~\tau^c>t).
$$
In section~\ref{sec-norm-propagators} we shall see that $\eta_t$ satisfies a nonlinear integro-differential equation given in weak form for any smooth functions by the formula
\begin{eqnarray}
\partial_t\eta_t(f)
&=&\eta_t(\La(f))-\eta_t(fV)+\eta_t(f)\eta_t(V)\label{evol-eta-cov}.
\end{eqnarray}

In contrast with conventional Markov processes, the flow of conditional probability measures
$\eta_t$ has a nonlinear evolution semigroup; that is, for any $s\leq t$ we have
\begin{equation}
\eta_t=\Phi_{t-s}\left(\eta_s\right),
\label{def-Phi}
\end{equation}
for some {\em nonlinear mapping} $\Phi_{t}$ from the set of probability measures on $\RR^r$ into itself. For a detailed description of these nonlinear transformations we refer to section~\ref{sec-norm-measures} and section~\ref{sec-norm-propagators}.

The stability analysis of the nonlinear evolution semigroups given by the composition of mappings $\Phi_{t+s}=\Phi_t\circ\Phi_s$ is closely related to the long time behavior of
a particle evolving in the ground state $h_0$, sometimes called the $h$-process, denoted $(X_t^h)_{t \ge 0}$ and defined by the
stochastic differential equation
\begin{equation}\label{def-X-h-d-intro}
dX^h_{t}=\left(AX^h_{t}+R\,\nabla \log h_0(X^h_{t})\right)~dt+BdW_t.
\end{equation}
The initial distribution of the random state $X^h_0$ is defined by a Boltzmann-Gibbs transformation of $\eta_0$ with respect to the ground state $h_0$; that is, we have that
\begin{equation}\label{def-BGT}
\eta^h_0= \BB_{h_0}(\eta_0)\quad \mbox{\rm with}\quad
 \BB_{h_0}(\eta_0)(dx):=\frac{1}{\eta_0(h_0)}~h_0(x)~\eta_0(dx),
\end{equation}
whenever $\eta_0(h_0)$ is a well-defined positive normalising constant.  By (\ref{def-h}), the generator $\La^h$  of the diffusion process $X_t^h$ is defined as $\La$ by replacing $A$ by the matrix $(A-RQ_{\infty})$.

The  distribution $\eta^h_t$ of the random states $X^h_t$ and  the distribution $\eta_t$ of a non-absorbed particle are connected
by the Boltzmann-Gibbs transformation; that is, for any time horizon $t\geq 0$ we have that
\begin{equation}\label{BG-h0-eta-intro}
\eta^h_t= \BB_{h_0}(\eta_t)\quad \mbox{\rm and}\quad \eta_t=\BB_{h_0^{-1}}(\eta_t^h)\quad \mbox{\rm with}\quad h_0^{-1}(x):=1/h_0(x).
\end{equation}
In the same vein, the Markov transitions of 
$X^h_t$ defined by the transition probabilities
$$
\Ka^h_t(x,dy):=\PP(X^h_t\in dy~|~X^h_0=x)
$$
are connected to the Feynman-Kac propagator $\Ka_t$ discussed in (\ref{def-Ka}) and (\ref{FKP-2}) by the formula
\begin{equation}\label{Ka-Kah} 
\exp{\left(\lambda_0 t\right)}~\Ka_t(x,dy)=~h_0(x)~\Ka^h_t(x,dy)~h_0^{-1}(y).
\end{equation}
The Boltzmann-Gibbs formulae (\ref{BG-h0-eta-intro}) and  (\ref{Ka-Kah}) remain valid for nonlinear diffusions $X_t$ (see for instance exercise 445 in~\cite{dmpenev}).
Under our controllability conditions (\ref{def-contr-obs}), we shall prove that the flow of probability measures $\eta_t$ and $\eta^h_t$ converge exponentially fast as $t\rightarrow\infty$  towards a pair of unique limiting measures, $\eta_{\infty}$ and $\eta^h_{\infty}$. That is for any $t\geq 0$ we have
\begin{equation}\label{cv-h0-proba-fixed point}
\eta_{\infty}=\Phi_{t}\left(\eta_{\infty}\right)\quad \mbox{\rm and}\quad \eta^h_{\infty}(dy)=(\eta^h_{\infty}\Ka^h_t)(dy):=\int \eta^h_{\infty}(dx)\Ka^h_t(x,dy).
\end{equation} 
The uniqueness property of $\eta_{\infty}$ is discussed at the end of section~\ref{gse-sec}.
The measure $\eta_{\infty}$ satisfies a nonlinear fixed point equation and it is sometimes called a quasi-invariant probability measure. 
Another important problem is to describe these limiting measures in terms of the parameters of the model and to quantify, with some precision, the convergence decays to equilibrium. 

 To briefly outline our answers to these questions, we denote by $\Na(m,P)$ an $r$-dimensional Gaussian probability measure with mean $m\in\RR^r$ and covariance matrix $P$.
 
 We also let  $(\widehat{X}_t,P_t)\in (\RR^r\times \RR^{r\times r})$ be the solution of the {\em coupled evolution equations} given by the system
\begin{equation}\label{ricc-intro}
\left\{\begin{array}{rcl}
 \partial_t\widehat{X}_t&=&(A-P_tS)\widehat{X}_t
\\
&&\\
\partial_tP_t&=&\mbox{\rm Ricc}(P_t)\quad \mbox{\rm with}\quad
 \mbox{\rm Ricc}(P):=AP + PA^{\prime}+R-PSP
\end{array}\right.\end{equation}
for some initial state $\widehat{X}_0\in\RR^r$ and some given positive semi-definite matrix $P_0$. 
Under the our controllability conditions (\ref{def-contr-obs}), we shall see that $\widehat{X}_t$ converges exponentially fast to $0$ as $t\rightarrow\infty$, and
the Riccati  matrix $P_t$ converges exponentially fast as $t\rightarrow\infty$ to a single positive fixed point matrix $P_{\infty}$ satisfying the continuous time algebraic Riccati equation $ \mbox{\rm Ricc}(P_{\infty})=0$.

Our main reason for introducing the coupled process $(\widehat{X}_t,P_t)$ comes from the following pivotal Gaussian preserving property
\begin{equation}\label{cv-eta-infty}
\eta_0=\Na(\widehat{X}_0,P_0)\Longrightarrow\forall t\geq 0\quad
 \eta_t=\Na(\widehat{X}_t,P_t)\longrightarrow_{t\rightarrow\infty}\eta_{\infty}:=\Na(0,P_{\infty})
\end{equation}

In addition, the zero-point energy level $\lambda_0$ of the Hamiltonian $\Ha$ is related to the survival probability of an non absorbed particle by the following formulae
\begin{equation}\label{survival-proba}
\begin{array}{l}
\displaystyle-\frac{1}{t}\log{\PP_{\eta_0}(\tau^c>t)}=\displaystyle\frac{1}{2t}\int_0^t\left(\widehat{X}_s^{\prime}\,S\,\widehat{X}_s+\tr(SP_s)\right)~ds. \\
\\
\displaystyle\longrightarrow_{t\rightarrow\infty}\lambda_0\stackrel{\forall s>0}{=}-\frac{1}{s}\log{\PP_{\eta_{\infty}}(\tau^c>s)}=\frac{1}{2}\,\tr(SP_{\infty})=\eta_{\infty}(V):=\int~\eta_{\infty}(dx) V(x).
\end{array}
\end{equation}

The Gaussian preserving property is discussed in Theorem~\ref{theo-eta-eta-c}.
Formula (\ref{survival-proba}) is a consequence of the exponential formula
 (\ref{gamma-1}) applied to the unit function.
The trace formula $\tr\left(R\,Q_{\infty}\right)=\tr(SP_{\infty})$ for non-necessarily reversible models  is proved in (\ref{key-QP-2}).
The convergence of the mean and covariance matrices $$(\widehat{X}_t,P_t)\rightarrow_{t\rightarrow\infty} (0,P_{\infty})
$$  can be made precise using the Lipschitz exponential decays estimates presented in section~\ref{ref-floquet-sec} and section~\ref{ref-sect-expo-cv}
 (see for instance Theorems~\ref{cor-phi-Lip} and~\ref{cor-X-hat-lip}).
The convergence of the distributions $$
 \eta_t^h\rightarrow_{t\rightarrow\infty} \eta^h_{\infty}\quad \mbox{\rm and}\quad
  \eta_t\rightarrow_{t\rightarrow\infty} \eta_{\infty}
 $$
can also be quantified with some precision in terms of relative entropy (cf. section~\ref{entropy-sec}) or in terms of Wasserstein distances (cf. section~\ref{wasserstein-sec}). See also the non-asymptotic estimates stated in theorem~\ref{theo-intro-entrop-wasserstein}.

Whenever the free particle $X_t$ is reversible with respect to some measure $\upsilon$ the quasi-invariant distribution  $\eta_{\infty}$ discussed above is related to the ground state $h_0$ by the Boltzmann-Gibbs formula
\begin{equation}\label{cv-h0-proba}
\eta_{\infty}=\BB_{h_0}(\upsilon)\quad \mbox{\rm and}\quad \eta^h_{\infty}= \BB_{h_0}(\eta_{\infty})=\BB_{h^2_0}(\upsilon).
\end{equation}
The proof of the above assertion  is provided in section~\ref{rev-section}. 
The left hand side assertion in (\ref{cv-h0-proba}) indicates that the stationary density of a non absorbed particle with respect to $\upsilon$  is proportional to  the ground state, while
the stationary distribution with respect to $\upsilon$ of a particle evolving in the ground state is proportional to the square of the ground state. These Boltzmann-Gibbs formulae yields the Hilbert space isometry
\begin{equation}\label{isometry}
\left\{
\begin{array}{l}
\Upsilon_h ~:~f\in \LL_2( \eta_{\infty}^h)\mapsto \Upsilon_h(f):= \sqrt{\eta_{\infty}^h(h^{-2}_0)}~f~h_0~\in \LL_2(\upsilon)\\
\\
 \mbox{\rm with the inverse}\quad\Upsilon_h^{-1}(f)=\sqrt{\upsilon(h^2_0)} ~f~h_0^{-1}.
 \end{array}\right.
 \end{equation}
This a direct consequence of the fact that \eqref{cv-h0-proba} implies that
$$
\langle f,g\rangle_{2,\upsilon}:=\upsilon(fg)=\upsilon(h_0^2)~\langle h_0^{-1}f,h_0^{-1}g\rangle_{2,\eta_{\infty}^h}\quad \mbox{\rm and}\quad \eta_{\infty}^h(h^{-2}_0)=1/\upsilon(h_0^2).
$$

\medskip

The Feynman-Kac propagator is connected to the semigroup of the 
particle evolving in the ground state $h_0$ via the operator formulae
$$
\exp{\left(\lambda_0 t\right)}~\Ka_t=\Upsilon_h\circ\Ka_t^h\circ\Upsilon_h^{-1}\quad
\mbox{\rm and}\quad -\La +(V-\lambda_0)=\Upsilon_h\circ\La^h\circ\Upsilon_h^{-1}.
$$ 
The Boltzmann-Gibbs mappings discussed in (\ref{BG-h0-eta-intro}), (\ref{Ka-Kah}) and (\ref{cv-h0-proba}) and the Hilbert space isometry (\ref{isometry}) allow to transfer directly any known regularity property at the level of the $h$-process $(\Ka^h_t,\eta^h_t,\eta^h_{\infty})$ to the Feynman-Kac model $(\Ka_t,\eta_t,\eta_{\infty})$, and vice versa.

\section{Basic notation and preliminary results}
\subsection{Some norms and matrix spaces}\label{prelim-matrix-sec}
We denote by $\Ma_{r_1,r_2}$ the set of $(r_1\times r_2)$-matrices with real entries and $r_1,r_2\geq 1$. When $r=r_1=r_2$ we write  $\Ma_{r}$ instead of $\Ma_{r,r}$ the set of square $(r\times r)$-matrices.

 A  square root  of a square matrix $A\in \Ma_{r}$ is a (non unique) matrix $A^{1/2}$ such that
$A^{1/2}A^{1/2}=A$. When $A$ has positive eigenvalues, we choose  the square root $A^{1/2}$ that has positive eigenvalues.
We let $\Sa_r\subset \Ma_{r}$ denote the subset of symmetric matrices, $\Sa_r^0\subset\Sa_r$ the subset of positive semi-definite matrices, and $\Sa_r^+\subset \Sa_r^0$ the subset of positive definite matrices. We also let $\Sa_r^-$ the set of negative definite matrices.

 Given $B\in \Sa_r^0-\Sa_r^+$ we denote by $B^{1/2}$ a (non-unique) but symmetric square root of $B$ (given by a Cholesky decomposition). When $B\in\Sa_r^+$ we choose the principal (unique) symmetric square root. We write $A^{\prime}$ to denote the transposition of a matrix $A$, and $A_{sym}=(A+A^{\prime})/2$ to denote the symmetric part of $A\in\Ma_{r}$. We denote by $\mbox{\rm Spec}(A)$ the spectrum of $A$ defined by
$$
\mbox{\rm Spec}(A):=\left\{\lambda~|~\lambda~\mbox{eigenvalue of}~A\right\},
$$
where each eigenvalue is listed the number of times it occurs as a root of the characteristic polynomial of $A$. We also denote by 
$\Ga l_r\subset \Ma_r$  the general linear group of invertible matrices.
The set $\Ma_{r}$ is equipped with the spectral or the Frobenius norms (a.k.a. Hilbert-Schmidt norm) defined by
 $$
 \Vert A \Vert=\sqrt{\lambda_{max}(AA^{\prime})}\leq \Vert A\Vert_{\tiny F}:=\sqrt{\tr(AA^{\prime})}\leq \sqrt{r}~\Vert A \Vert,
$$ 
where $\lambda_{max}(\cdot)$ denotes the maximal eigenvalue. The minimal eigenvalue is denoted by $\lambda_{min }(\cdot)$. Let $\tr(A)=\sum_{1\leq i\leq r}A(i,i)$ denote the trace operator. We also denote by $\mu(A)=\lambda_{max}(A_{sym})$ the logarithmic norm and by 
$$
\varsigma(A):=\max_{\lambda\in\mbox{\rm Spec}(A)}{\left\{\mbox{\rm Re}(\lambda)\right\}},
$$ the spectral abscissa. 
We recall that 
$$
\Vert A\Vert\geq \mu(A)=\lambda_{max}(A_{sym})\geq \varsigma(A):=\max{\left\{\mbox{\rm Re}(\lambda)~:~\lambda\in \mbox{\rm Spec}(A)\right\}}.
$$
 A matrix $A$ is said to be stable (a.k.a. Hurwitz) when $\varsigma(A)<0$. We recall that
 \begin{equation}\label{stable-matrix-exp}
\varsigma(A)<0\Longrightarrow \exists \alpha,\beta>0~:~\forall t\geq 0\quad
 \Vert e^{tA}\Vert \,\leq\, \alpha\,e^{-\beta t}.
\end{equation}
 The parameters $(\alpha,\beta)$ can be made explicit in terms of the spectrum of the matrix $A$. For instance, applying Coppel's inequality (cf. Proposition 3 in~\cite{coppel1978stability}), for any $t\geq 0$ and any  for any $0<\gamma<1$ we can choose 
 $$
 \alpha=({a}/{\gamma})^{r-1}\quad \mbox{\rm and} \quad\beta=(1-\gamma)\varsigma(A)\quad \mbox{\rm with} \quad a:=2{\Vert A\Vert}/{\vert\varsigma(A)\vert}.
$$ 
In the case where $\beta=-\mu(A)>0$, we can choose $\alpha=1$.

\subsection{Relative entropy and metrics}
The $n$-th Wasserstein distance between two probability measures
$\nu_1$ and $\nu_2$ on $\RR^r$ is defined for any parameter $n\geq 1$ by the formula
$$
\WW_n(\nu_1,\nu_2)=\inf{\left\{\EE\left(\Vert Z_1-Z_2\Vert^n\right)^{\frac{1}{n}}
\right\}}.
$$
The infimum in the above display is taken over all pairs of random variables $(Z_1,Z_2)$
such that $\mbox{\rm Law}(Z_i)=\nu_i$, for $i=1,2$. 

We denote by $ \mbox{\rm Ent}\left(\nu_1~|~\nu_2\right)$
the Boltzmann-relative entropy, defined as
$$
 \mbox{\rm Ent}\left(\nu_1~|~\nu_2\right):=\int~\log\left(\frac{d\nu_1}{d\nu_2}\right)~d\nu_1,
 $$
whenever $\nu_1\ll\nu_2$, and $+\infty$ otherwise.
Further, the Fisher information is defined by
$$
 \Ja\left(\nu_1~|~\nu_2\right):=\int~\Vert\nabla \log\left(\frac{d\nu_1}{d\nu_2}\right)\Vert^2~d\nu_1,
 $$
 if $\log{d\nu_1/d\nu_2}\in \LL_2(\nu_1)$, and $+\infty$ otherwise.
The total variation distance between the measures $\nu_1$ and $\nu_2$ is defined by
$$
\Vert \nu_1-\nu_2\Vert_{\tiny tv}:=\frac{1}{2}\sup{\left\{\vert\nu_1(f)-\nu_2(f)\vert~:~f~\mbox{\rm s.t.}~\Vert f\Vert_{\infty}\leq 1\right\}},
$$
with the uniform norm and Lebesgue integrals defined, respectively, by
$$
\Vert f\Vert_{\infty}:=\sup_{x\in\RR^r}\vert f(x)\vert, \qquad 
\nu_i(f):=\int~\nu_i(dx)~f(x).
$$
Finally, given some locally finite measure $\nu$ on $\RR^r$, for any $n\geq 1$ 
we denote by $\LL_n(\nu)$ the Banach space of measurable functions $f$ on $\RR^r$
equipped with the norm
$$
\Vert f\Vert_{n,\nu}:=\nu\left(\vert f\vert^n\right)^{1/n}.
$$

\subsection{Evolution semigroups}
The evolution semigroup $(\widehat{X}_t(x,P),\phi_t(P))$ of the equations (\ref{ricc-intro}) starting at $(x,P)$ satisfies the coupled equations
\begin{equation}\label{coupled-equations}
\left\{
\begin{array}{rclclcl}
 \partial_t\widehat{X}_t(x,P)&=&(A-\phi_t(P)S)\widehat{X}_t(x,P)& \mbox{\rm with}&\widehat{X}_0(x,P)&=&x\\
&&&&&&\\
\partial_t\phi_t(P)&=&\mbox{\rm Ricc}(\phi_t(P))& \mbox{\rm with}&\phi_0(P)&=&P.
\end{array}\right.
\end{equation}
By the Gaussian preserving property of the measure-valued nonlinear semigroup $\Phi_t$, defined in \eqref{def-Phi}, we have
\begin{equation}\label{Phi-Na}
\Phi_t\left(\Na(x,P)\right)=\Na\left(\widehat{X}_t(x,P),\phi_t(P)\right)\quad \mbox{\rm and}\quad
\Phi_t(\delta_x)=\Na\left(\widehat{X}_t(x,0),\phi_t(0)\right).
\end{equation}
We refer to Theorem~\ref{theo-eta-eta-c}, as well as section~\ref{ref-floquet-sec} and Proposition~\ref{alternative-description-xw-P} for different ways of writing the evolution semigroup $\left(\widehat{X}_t(x,P),\phi_t(P)\right)$.

The controllability conditions (\ref{def-contr-obs}) are well known in filtering theory, see for instance \cite{Bittanti91,Lancaster1995} and the more recent articles \cite{ap-2016,bp-21}, and references therein. They ensure the existence of an unique pair $(P^-_{\infty},P_{\infty})$
of negative and positive fixed point matrices of the algebraic Riccati equation
\begin{equation}\label{intro-ricc}
  \mbox{\rm Ricc}(P^-_{\infty})=0=  \mbox{\rm Ricc}(P_{\infty})\end{equation}
In addition, the matrices
\begin{equation}\label{CARE}
A-P_{\infty}S~~\mbox{\rm and}~~ A^{\prime}+(P^-_{\infty})^{-1}R
\end{equation}	
are stable (a.k.a. Hurwitz).

For a more thorough discussion on the above assertions we refer to~\cite[Chapter 3]{Bittanti91},~\cite{Lancaster1995} and the more recent articles \cite{ap-2016,bp-21}. 

As already mentioned, the pair of matrices $((P^-_{\infty})^{-1},P_{\infty}^{-1})$ satisfy the same fixed point equation as tthat of $(P^-_{\infty},P_{\infty})$ by replacing $(A,R,S)$ by $(-A^{\prime},S,R)$. In addition, the matrices $(Q_{\infty}^-,Q_{\infty})$ defined by
 \begin{equation} \label{formula-example-Pinfty-0-Q}
Q_{\infty}^-:=-P^{-1}_{\infty}< 0<
 Q_{\infty}:=-(P^-_{\infty})^{-1}
\end{equation}
satisfy the same fixed point equation as $((P^-_{\infty})^{-1},P_{\infty}^{-1})$  by replacing $A$ by $(-A)$. Thus, the matrices $(Q_{\infty}^-,Q_{\infty})$ satisfy the algebraic Riccati matrix equation (\ref{def-h}).

Let $\Ea_{s,t}(P)$ be the exponential semigroup associated with the matrix flow $u\mapsto (A-\phi_{u}(P)S)$; that is
the solution for any $0\leq s\leq t$ of the matrix evolution equations
\begin{equation}\label{syst-E-def}
 \partial_t \Ea_{s,t}(P)=(A-\phi_t(P)S)\,\Ea_{s,t}(P)\quad\mbox{\rm and}\quad
\partial_s \Ea_{s,t}(P)=-\Ea_{s,t}(P)\,(A-\phi_s(P)S),
\end{equation}
with $\Ea_{s,s}(P)=I$ and where we often write $\Ea_{t}(P)$ for $\Ea_{0,t}(P)$. In this notation, the solution of the right hand side equation in (\ref{ricc-intro}) takes the form
\begin{equation}\label{Est}
\widehat{X}_t(x,P)=\Ea_t(P)x\quad \mbox{\rm and}\quad
\Ea_{s,s+t}(P_{\infty})=\Ea_{t}(P_{\infty})=\exp{\left(t(A-P_{\infty}S)\right)}.
\end{equation}

 From a mathematical viewpoint, it is tempting to integrate sequentially the differential equations (\ref{syst-E-def}), to obtain an explicit description of $\Ea_{s,t}(P)$ in terms of the Peano-Baker series~\cite{Peano,baker1905}, see also~\cite{Brockett,Frazer,Ince}. Another natural strategy is to express the semigroup as a true
 matrix exponential involving a Magnus series expansion of iterated integrals on the Lie algebra generated by the matrices $(A-\phi_u(P)S)$, with $s\leq u\leq t$. For more details on these exponential expansions we refer to~\cite{blanes,magnus}. In practical terms, the use of Peano-Baker and/or exponential Magnus series in the study of the stability properties of time-varying linear dynamical systems is rather limited.

For any semi-definite positive initial matrix $P\in \Sa_r^0$, we have the following results (see e.g. \cite{ap-2016,bucy2}),
\begin{equation}\label{exp-stab-phi}
\begin{array}{l}
\forall t\geq \delta >0, \quad 0<\Pi_{-,\delta}\leq \phi_t(P)\leq \Pi_{+,\delta}
\quad\mbox{\rm and}\quad
 \forall t\geq 0,\quad
 \Vert \Ea_{t}(P_{\infty}) \Vert\,\leq\, \alpha\,e^{-\beta t},
\end{array}\end{equation}
for some positive matrices $\Pi_{-,\delta},\Pi_{+,\delta}$ and some $\alpha,\beta>0$, all of which depend on the model parameters $(A,R,S)$. The right hand side assertion comes from the fact that $(A-P_{\infty}S)$ is a stable matrix. Here, $\Vert \cdot \Vert$ stands for the spectral norm of matrices.

\subsection{Linear diffusion processes}
Consider an $r$-dimensional process given by the linear stochastic differential equation
\begin{equation}\label{lin-Gaussian-diffusion-filtering}
dX_t=AX_tdt+BdW_t
\end{equation}
for some initial state $X_0$ with distribution $\eta_0$ on $\RR^r$.
In the above display, $W_t$ is an $r_1$-dimensional Brownian motion, for some $r_1\geq 1$, which is independent of $X_0$, and $B$ is a $(r\times r_1)$-matrix such that $BB^{\prime}=R$. The random state vectors $X_t$ and the Brownian states $W_t$ are column vectors. 
 The generator of the diffusion process  $X_t$ coincides with  the second order differential kinetic energy operator defined in (\ref{def-La-intro}). In the same vein, 
the $h$-process (\ref{def-X-h-d-intro}) satisfies the
stochastic differential equation
\begin{equation}\label{def-X-h-d}
dX^h_{t}=\left(A-RQ_{\infty}\right)X^h_{t}~dt+BdW_t.
\end{equation}
Due to (\ref{CARE}), the matrix $(A-RQ_{\infty})$ is Hurwitz so that
$X^h_{t}$ is a stable Ornstein-Uhlenbeck process even when $A$ is unstable. This property ensures the existence of some parameters $\alpha_h,\beta_h>0$ such that
\begin{equation}\label{expo-decay-h-exp}
\Vert  e^{(A-RQ_{\infty})t}\Vert\leq \alpha_h~e^{-\beta_h t}\quad \mbox{\rm and}\quad
\iota_h:=\int_0^{\infty}\Vert  e^{(A-RQ_{\infty})t}\Vert^2~dt\leq \frac{\alpha_h^2}{2\beta_h}.
\end{equation}
For a more thorough discussion on the exponential decays of matrix exponential-type semigroups (a.k.a. fundamental matrices) we refer to section~\ref{prelim-matrix-sec}. Recall that
$$
\eta^h_0=\Na(x,P)\Longrightarrow\forall t\geq 0\quad \eta^h_t=\Na(\widehat{X}^h_t(x),\phi^h_t(P)),
$$  
with the mean vector
\begin{equation}\label{Xh-mean-cov}
\widehat{X}^h_t(x)=e^{(A-RQ_{\infty})t}x,
\end{equation}
and the covariance matrices
$$
\phi^h_t(P):=e^{(A-RQ_{\infty})t}Pe^{(A-RQ_{\infty})^{\prime}t}+\int_0^t
e^{(A-RQ_{\infty})s}Re^{(A-RQ_{\infty})^{\prime}s}~ds.
$$
Notice that $\phi_t^h(P)$ is the evolution semigroup of associated with the matrix valued differential equation 
$$
\partial_tP^h_t=(A-RQ_{\infty})P^h_t+P^h_t(A-RQ_{\infty})^{\prime}+R\quad \mbox{\rm with}\quad P^h_0=P.
$$
The stochastic flow $X^h_t(x)$ of the $h$-process is defined as in (\ref{def-X-h-d}) by choosing  the initial condition $X^h_0(x)=x$. 
The random function
$x\mapsto X^h_t(x)$  can be seen as the Gaussian random field 
$$
X^h_t(x)\sim \Na(\widehat{X}^h_t(x),\phi^h_t(0)).
$$

Without additional conditions on the matrices $(A,R)$  the $h$-process defined in (\ref{def-X-h-d}) may not be reversible, see for instance Theorem~\ref{ref-theo-h}.
Nevertheless, for any $x\in\RR^d$ and $P\geq 0$, we have
  \begin{equation}\label{def-Ph}
\Na(\widehat{X}^h_t(x),\phi^h_t(P))
 \longrightarrow_{t\rightarrow \infty}  \eta^h_{\infty}:=\Na(0,P^h_{\infty}),
\end{equation}
with the limiting covariance matrix
$$
P^h_{\infty}:=
\int_0^{\infty}
e^{(A-RQ_{\infty})s}R\,e^{(A-RQ_{\infty})^{\prime}s}~ds=(P_{\infty}^{-1}+Q_{\infty})^{-1}.
$$
The right hand side assertion is a consequence of the Gramian formula (\ref{Delta-P-Q}) (see also Theorem~\ref{theo-eta-eta-c}).
Exponential decay estimates to equilibrium can be easily extracted from the exponential inequalities (\ref{expo-decay-h-exp}).

{\color{black}\subsection{Reversible models}\label{rev-models-intro}
Without further mention, until the end of this section we shall assume that matrices  $(A,R,S)$ satisfying the rank condition (\ref{def-contr-obs}). In addition, 
we have  $R>0$ and $AR=RA^{\prime}$.

In this situation, the diffusion  process $X_t$ defined  in (\ref{lin-Gaussian-diffusion-filtering}) is reversible with respect to the locally finite measure
\begin{equation}\label{def-upsilon-intro}
\upsilon(dx):=\exp{\left(U(x)\right)}~dx\quad \mbox{\rm with}\quad U(x):=x^{\prime}R^{-1}Ax.
\end{equation}
The stochastic differential equations (\ref{lin-Gaussian-diffusion-filtering}) and (\ref{def-X-h-d}) resume to the Langevin diffusions
\begin{equation}\label{h-process-reversible}
dX_t=\frac{1}{2}\,R\,\nabla U(X_t) dt+B dW_t\quad\mbox{\rm and}\quad
dX^h_{t}=\frac{1}{2}~R\,\nabla\left(U+ \log h_0^2\right)(X^h_{t})~dt+BdW_t.
\end{equation}
All the limiting covariance matrices $(P_{\infty},P_{\infty}^h,Q_{\infty})$ can be explicitly computed in terms of the parameters $(A,R,S)$. The matrix $Q_{\infty}$ was already given in (\ref{ref-Weyls-2-intro-0}) and we have
$$
\frac{1}{2}~(P^h_{\infty})^{-1}=
Q_{\infty}-R^{-1}A\quad \mbox{\rm and}\quad P_{\infty}^{-1}=Q_{\infty} - 2 R^{-1}A.
$$
For a more thorough discussion on these formulae we  refer to
section~\ref{rev-section}.

In the reversible case, it is convenient to rewrite the generator of the $h$-process given by (\ref{h-process-reversible}) in the divergence form
$$
\La^h(f)=\frac{1}{2}~e^{U_h}~\sum_{1\leq i\leq r}~\partial_{x_i}\left(e^{-U_h}\partial_{x_i}f\right)\quad \mbox{\rm with}\quad U_h(x):=\frac{1}{2}x^{\prime}(P^h_{\infty})^{-1}x.
$$

\begin{theo}\label{ref-theo-h}
For any $t\geq 0$ we have the master equation
\begin{equation}\label{intro-master}
\eta^h_{\infty}(dx)~\Ka_t^h(x,dy) = \eta^h_{\infty}(dy)~\Ka_t^h(y,dx)
\end{equation}
with the distribution $\eta^h_{\infty}$ defined in (\ref{def-Ph}). In addition, we have the density-transport formulae
\begin{equation}\label{density-transport}
\eta_0^h(dx):=f_0(x)~\eta^h_{\infty}(dx)\Longrightarrow \eta_t^h(dx)=f_t(x)~\eta^h_{\infty}(dx)\quad \mbox{with}\quad f_t(x):=\Ka_t^h(f_0)(x).
\end{equation}

\end{theo}
\proof
 Using the divergence form of the generator we check that
\begin{equation}\label{ipb-ref}
\eta^h_{\infty}\left(g~\La^h(\Ka_t^h(f_0)) \right)=-\frac{1}{2}~~\int~\eta^h_{\infty}(dx)~\sum_{1\leq i\leq r}\partial_{x_i}(g)(x)~\partial_{x_i}\Ka_t^h(f_0)(x)
\end{equation}
for sufficiently smooth functions $f$ and $g$ for which we can perform integration by parts.
This yields for any $f,g\in \LL_2(\eta^h_{\infty})$ the formula
$$
\eta^h_{\infty}(f\,\Ka_t^h(g)) = \eta^h_{\infty}(\Ka_t^h(f)\, g),
$$
which is equivalent to \eqref{intro-master}. The density-transport formulae (\ref{density-transport}) is a direct consequence of the reversible property (\ref{intro-master}).
\cqfd

The convergence to equilibrium of Langevin-type $h$-processes can be studied in terms of the Boltzmann-relative entropy using the rather well-known de Bruijn identity
\begin{equation}\label{de-Bruijn}
\partial_t \mbox{\rm Ent}\left(\eta_t^h~|~\eta_{\infty}^h\right)=-\frac{1}{2}~\Ja\left(\eta_t^h~|~\eta_{\infty}^h\right).
\end{equation}
From this, one can obtain the exponential decays of the Fisher information
\begin{equation}\label{exp-decay-Fisher}
\Ja\left(\eta_t^h~|~\eta_{\infty}^h\right)\leq \Vert  e^{(A-RQ_{\infty})t}\Vert^2~\Ja\left(\eta_0^h~|~\eta_{\infty}^h\right),
\end{equation}
which also yields the log-Sobolev inequality
\begin{equation}\label{log-Sobolev}
\mbox{\rm Ent}\left(\eta_0^h~|~\eta_{\infty}^h\right)\leq \frac{\iota_h}{2}~\Ja\left(\eta_0^h~|~\eta_{\infty}^h\right),
\end{equation}
where the parameter $\iota_h$ introduced in (\ref{expo-decay-h-exp}).
Applying the log-Sobolev inequality to $\eta^h_t$, the
 de Bruijn identity now yields the free energy exponential decays
\begin{equation}\label{free-energy-decays}
\partial_t \mbox{\rm Ent}\left(\eta_t^h~|~\eta_{\infty}^h\right)=-\frac{1}{2}~\Ja\left(\eta_t^h~|~\eta_{\infty}^h\right)
\leq -\iota_h^{-1}~\mbox{\rm Ent}\left(\eta_t^h~|~\eta_{\infty}^h\right).
\end{equation}
The proofs of the assertions \eqref{de-Bruijn}-\eqref{free-energy-decays} follow standard probabilistic manipulations and are thus provided in the Appendix. 
We summarise the above discussion with the following theorem.
\begin{theo}\label{theo-entropy-fisher}
For any $t\geq 0$ we have the relative entropy exponential decays
$$
\Ja\left(\eta_t^h~|~\eta_{\infty}^h\right)\leq  \alpha^2_h~e^{-2\beta_h t}\Ja\left(\eta_0^h~|~\eta_{\infty}^h\right)\quad \mbox{and}\quad 
\mbox{\rm Ent}\left(\eta_t^h~|~\eta_{\infty}^h\right)\leq e^{-t/\iota_h}~\mbox{\rm Ent}\left(\eta_0^h~|~\eta_{\infty}^h\right)
$$
with the parameters $(\iota_h,\alpha_h,\beta_h)$ introduced in (\ref{expo-decay-h-exp}).
\end{theo}
From the practical point of view, the functional inequalities discussed above are rarely useful when the matrix $Q_{\infty}$ and thus the limiting measure $\eta_{\infty}^h$ is not explicitly known.
}

\section{Statement of some main results}\label{sec-main-results}
\subsection{Ground state energy}\label{gse-sec}
Our first main result provides an explicit description of the ground state energy of the Hamiltonian operator for general matrices $(A,R,S)$ in terms of the negative and positive fixed points of the algebraic Riccati equation (\ref{intro-ricc}). 
\begin{theo}\label{theo-1-intro}
For any matrices  $(A,R,S)$ satisfying the rank condition (\ref{def-contr-obs}), 
the function $h_0$ in (\ref{def-h}) is the ground state of the Hamiltonian  $\Ha$ introduced in  (\ref{def-H-intro}); that is, we have 
\begin{equation}\label{intro-h0}
\Ha(h_0)=\lambda_0\, h_0\quad \mbox{with }\quad \lambda_0:=\frac{1}{2}\,\tr\left(S P_{\infty}\right)=\frac{1}{2}\,\tr\left(R\,Q_{\infty}\right)>0.
\end{equation}
In addition, we have the Feynman-Kac propagator formula stated in (\ref{Ka-Kah}).
\end{theo}

The proof of the right hand side trace formula in (\ref{intro-h0}) is provided in section~\ref{section-Riccati} dedicated to Riccati algebraic equations (see Lemma~\ref{lem-QP-2}). The first assertion and the Feynman-Kac propagator formula stated in (\ref{Ka-Kah}) is a direct consequence of (\ref{def-Ka}) and Theorem~\ref{theo-c-h0}.

\begin{theo}\label{theo-eta-eta-c}
For any matrices $(A,R,S)$, the law $\eta^h_t$ of the random states $X^h_t$ of the $h$-process defined in (\ref{def-X-h-d}) and the distribution of the non-absorbed particle defined in (\ref{def-d-non-abs}) are connected by the Boltzmann-Gibbs transformation (\ref{BG-h0-eta-intro}).
In addition, we have the Gaussian preserving property
$$
\eta_0=\Na(\widehat{X}_0,P_0)\quad\Longrightarrow\quad \forall t\geq 0\qquad
\eta_t^h=\Na(\widehat{X}_t^h,P_t^h)\quad \mbox{and}\quad \eta_t=\Na(\widehat{X}_t,P_t)
$$
with the parameters $(\widehat{X}_t,P_t)$ defined in (\ref{ricc-intro}), the covariance matrix
\begin{equation}\label{ref-link-h-wX}
P^h_t=(P_t^{-1}+Q_{\infty})^{-1}\quad\mbox{and the mean vector}\quad
\widehat{X}^h_t=P^h_t P_t^{-1}\widehat{X}_t.
\end{equation}
\end{theo}

The proof of Theorem~\ref{theo-eta-eta-c} is provided in the end of section~\ref{particle-absorption-sec}. The next theorem provides a non asymptotic expansion of the Feynman-Kac propagator.

\begin{theo}\label{theo-lamda0-K}
 For any time horizon $t\geq \delta>0$ and any $f\in \LL_1(\eta_{\infty})$ we have $$
e^{\lambda_0t}~\Ka_t(f)(x)
=\frac{h_0(x)}{\eta_{\infty}(h_0)}~\left(\eta_{\infty}(f)+\epsilon_t(f)(x)\right)~k_t(x)
$$
where 
$$
\epsilon_t(f)(x):=\Phi_t(\delta_x)(f)-\eta_{\infty}(f)
$$
and $k_t$ is a function satisfying
\begin{equation}\label{estim-final}
\exp{\left(-c_{\delta}\, e^{-2\beta t}
\right)}\leq k_t(x)\leq \exp{\left(c_{\delta} ~(1+\Vert x\Vert^2)~e^{-2\beta t}\right).
}
\end{equation}
In the above display  $\beta$ stands for the parameter defined in (\ref{exp-stab-phi}), and $c_{\delta}<\infty$ is some finite constant whose value only depends on $\delta$.
\end{theo}
The proof of the above theorem is provided in section~\ref{theo-lamda0-K-proof}.

The convergence of $\eta_t$ to the limiting measure $\eta_{\infty}$ discussed in \eqref{cv-eta-infty} can be studied in terms of both the stationary properties of the $h$-process and the stability properties of the Riccati matrix flow $P_t$. Due to the exponential semigroup formula \eqref{Est}, the long time behavior of the mean vector $\widehat{X}_t$ is also directly related to $P_t$. Thus, in section~\ref{section-Riccati} we provide a brief discussion on Riccati matrix flows, including the Floquet-type theory developed in~\cite{bp-21}, as well as several Lipschitz type inequalities and exponential type decays to equilibrium for Riccati flows and their associated exponential semigroups. Applied to our context, these quantitative estimates allow one to prove a variety of non asymptotic convergence theorems. 

To give a flavour of these results, consider the initial distributions
\begin{equation}\label{spec-initial-Gauss}
\eta_0=\Na(x,P)\quad \mbox{and}\quad\mu_0=\Na(y,Q),
\end{equation}
for some $x,y\in\RR^r$ and some covariance matrices $P,Q\in\Sa^0_r$. Our main results can be summarised with the following theorem.
\begin{theo}\label{theo-intro-entrop-wasserstein}
There exists $\d = \d(P, Q) > 0$, which depends on the distance $P - Q$, such that for any $t\geq \delta$ we have
$$
\WW_2\left(\Phi_t(\eta_0),\Phi_t(\mu_0)\right)\leq c_{\delta}~e^{-\beta t}~\left(
\Vert x-y\Vert+(\Vert x\Vert\vee 1)~\Vert P-Q\Vert \right)
$$
and for sufficiently large time horizon we have
$$
\begin{array}{l}
\displaystyle
 \mbox{\rm Ent}\left(\Phi_t(\eta_0)~|~\Phi_t(\mu_0)\right)\leq c_{\delta}\left(~\Vert P-Q\Vert+~ \left(\Vert x\Vert^2~\Vert P-Q\Vert^2+\Vert x-y\Vert^2
\right)\right)~e^{-2\beta t}
\end{array}$$
for some finite constant $c_{\delta}$ and the parameter $\beta>0$ introduced in \eqref{exp-stab-phi}. 

In addition, when $P=Q$ for any $n\geq 1$ we have
$$
\WW_n\left(\Phi_t(\eta_0),\Phi_t(\mu_0)\right) \leq c_{\delta}~e^{-\beta t}~\Vert x-y\Vert
\quad
\mbox{and}
\quad
\mbox{\rm Ent}\left(\Phi_t(\eta_0)~|~\Phi_t(\mu_0)\right)\leq c_{\delta}~e^{-2\beta t}~\Vert x-y\Vert^2.
$$
\end{theo}
The case $P=Q$ is a direct consequence of the Gaussian preserving property (\ref{Phi-Na}) and the Lipschitz estimates stated in Theorem~\ref{cor-X-hat-lip}. For the general case, we refer the reader to Theorems~\ref{theo-time-delta-n} and~\ref{ref-W-intro}, where a more precise description of the constant $c_{\delta}$ and the time horizon in the relative entropy estimates are provided.

Total variation estimates for initial Gaussian measures can be deduced directly from the relative entropy estimates stated in theorem~\ref{theo-intro-entrop-wasserstein} using Pinsker's inequality 
$$
\Vert \Phi_t(\eta_0)-\Phi_t(\mu_0)\Vert_{\tiny tv}\leq\sqrt{\frac{1}{2}\,\mbox{\rm Ent}\left(\Phi_t(\eta_0)~|~\Phi_t(\mu_0)\right)}
$$
More generally (cf. Theorem~\ref{theo-non-Gauss}), for any pair of probability measures $\eta_0$ and $\mu_0$ on $\RR^r$ and any time horizon $t\geq \delta>0$ we have
$$
 \Vert\Phi_t(\eta_0)-\Phi_t(\mu_0)\Vert_{\tiny tv}\leq c_{\delta}(\eta_0,\mu_0)~e^{-\beta t}m
$$
where the parameter $\beta$ was introduced in (\ref{exp-stab-phi}), and
$c_{\delta}(\eta_0,\mu_0)$ is a finite constant that depends on the parameters $(\delta,\eta_0,\mu_0)$. The above result implies the uniqueness of the fixed point Gaussian distribution $\eta_{\infty}$ introduced in (\ref{cv-eta-infty}). Choosing $\mu_0=\eta_{\infty}$, it also shows that for any initial distribution $\eta_0$ the probability measure $\Phi_t(\eta_0)$ converges exponentially fast towards a Gaussian distribution as the time horizon $t\rightarrow\infty$. 

Theorem~\ref{theo-intro-entrop-wasserstein} also provides several ways of estimating the difference $\epsilon_t(f)(x)$ defined in Theorem \ref{theo-lamda0-K}. For instance, for any Lipschitz function $f$ with unit Lipschitz constant or for any bounded functions $g$ with unit uniform norm, we have the estimates
$$
\vert \epsilon_t(f)(x)\vert\leq \WW_1\left(\Phi_t(\delta_x),\eta_{\infty}\right)
\quad \mbox{\rm and}\quad
\vert\epsilon_t(g)(x)\vert\leq 2~\Vert\Phi_t(\delta_x)-\eta_{\infty}\Vert_{\tiny tv}.
$$

\subsection{Spectral theorems}\label{spectral-intro-sec}
Assume that matrices  $(A,R,S)$ satisfy (\ref{def-contr-obs}) and we have  $R>0$ and $AR=RA^{\prime}$.
Let $\Lambda^h$ be the matrix defined by
$$
\Lambda^h:=-(P^{h}_{\infty})^{-1/2}(A^2+RS)^{1/2}(P^{h}_{\infty})^{1/2},
$$
where the positive matrix $P^{h}_{\infty}>0$ was introduced in (\ref{def-Ph}).
We denote by $\Za:=\left(z_1,\ldots,z_r\right)$ the orthogonal matrix with columns given by the orthonormal eigenvector $z_i$ of the matrix $\Lambda^h$ associated with an eigenvalue $\lambda_i(\Lambda^h):=-\lambda_i^h<0$, for $i\in\{1,\ldots,r\}$. 

Under the reversibility condition, the matrix $(A^2+RS)$ {\em may not be symmetric} but it has positive eigenvalues
(see for instance  (\ref{ref-positive-A2}) and (\ref{spec-Lambda-h})), so that
the parameters $\lambda^h_i$ coincide with the square roots of the eigenvalues of  $(A^2+RS)$. 
We shall assume that these eigenvalues  are ranked in increasing order
 $$
 \lambda^h_1\leq \ldots\leq \lambda^h_r.
 $$

Let $\HH^r_n(x)$ be the collection of multivariate Chebychev-Hermite polynomials on $\RR^r$
indexed by the multiple indices $n=(n_1,\ldots,n_r)\in \NN^r$. We use the multiple index notation $n!=n_1!\times\ldots\times n_r!$. Denote by $(0)$ the multiple index with null entries $n_i=0$,
and recall that $\HH_{(0)}^r(x)=1$. Further, set $\lambda_0(\Lambda^h)=-\lambda_0^h:=0.$

We are now in position to state the main result of this section.
\begin{theo}\label{spectral-theo-intro}
 For any time horizon $t\geq 0$, we have the $\LL_2(\eta^h_{\infty})$-spectral decomposition
$$
\Ka^h_t(x,dy)=\sum_{n\in\NN^r}~e^{-\lambda_n^ht}~\varphi^h_n(x)\varphi^h_n(y)~\eta^h_{\infty}(dy)
$$
with the $\LL_2(\eta_{\infty}^h)$ orthonormal basis eigenfunctions $\varphi^h_n$ and corresponding eigenvalues $\lambda_n^h$ given respectively by
\begin{equation}\label{intro-def-varphi-ln}
\varphi^h_n(x):=\frac{1}{\sqrt{n!}}~
\HH^r_n\left(\Za^{\prime}(P^{h}_{\infty})^{-1/2}x\right),
\qquad\lambda_n^h:=\sum_{1\leq i\leq r}n_i~\lambda_i^h.
\end{equation}
\end{theo}
The proof of the above theorem is provided in section~\ref{sec-spectral-proofs}.

\medskip 

Thus, for any $n\in \NN^r$ and any $t\geq 0$ and $x\in\RR^r$ we have the formulae
\begin{equation}\label{Kh-varphi}
\Ka^h_t\left(\varphi^h_n\right)=e^{-\lambda_n^ht}~\varphi^h_n\quad \mbox{\rm and}\quad
\La^h\left(\varphi^h_n\right)=-\lambda_n^h\,\varphi^h_n.
\end{equation}
Exponential decays to equilibrium can be extracted directly from the spectral decomposition. For instance we have the following estimates.
\begin{cor}\label{Poincare-cor}
For any time horizon $t\geq 0$ we have the exponential decays to equilibrium
\begin{equation}\label{expo-decays}
\Vert \Ka^h_t(f)-\eta_{\infty}^h(f)\Vert_{2,\eta_{\infty}^h}\leq e^{-\lambda_1^ht}~\Vert f-\eta_{\infty}^h(f)\Vert_{2,\eta_{\infty}^h}.
\end{equation}
Equivalently, we have the Poincar\'e inequality 
\begin{equation}\label{poincare-ineq}
\lambda_1^h~\Vert f-\eta_{\infty}^h(f)\Vert_{2,\eta_{\infty}^h}^2\leq E_h(f,f):=-\eta_{\infty}^h(f\,\La^h(f)).
\end{equation}
\end{cor}
The proof of the above corollary is also provided in section~\ref{sec-spectral-proofs}. 

\medskip

The isometry (\ref{isometry}) shows that the $\LL_2(\upsilon)$ orthonormal basis and those of $\LL_2(\eta_{\infty}^h)$ are linked to each other by the formulae
$$
\varphi_n=\Upsilon_h \left(\varphi^h_n\right)\quad \mbox{\rm and}\quad 
\varphi^h_n=\Upsilon_h^{-1}\left(\varphi_n\right), \qquad n\in\NN_r.
 $$
{\color{black}The master equation (\ref{intro-master}) is equivalent to the reversible property 
\begin{equation}\label{rev-Kt}
\upsilon(dx)~\Ka_t(x,dy) =\upsilon(dy)~\Ka_t(y,dx)
\end{equation}
}

Rewritten in terms of Feynman-Kac propagators Theorem~\ref{spectral-theo-intro} takes the following form.
\begin{theo}\label{theo-Ka-spec}
For any $t\geq 0$ we have the $\LL_2(\upsilon)$ spectral decomposition
\begin{equation}\label{spectral-Ka-ref-app}
\Ka_{t}(x,dy)=\sum_{n\in\NN^r} e^{-\lambda_nt}~\varphi_n(x)~\varphi_n(y)~\upsilon(dy),
\end{equation}
with the $\LL_2(\upsilon)$ orthonormal basis  given for any $n\in\NN^r$ by the eigenstates
 $$
\varphi_n(x)=\left(\frac{1}{(2\pi)^r}~\mbox{det}\left(P_{\infty}^{-1}+Q_{\infty}\right)\right)^{1/4}
~h_0(x)\,\varphi^h_n(x),
 $$
 and corresponding eigenvalues $\lambda_n:=\lambda_0+\lambda^h_n$.
\end{theo}

Using the above spectral decomposition we check the formulae
$$
\Ka_t\left(\varphi_n\right)=e^{-\lambda_nt}~\varphi_n\quad \mbox{\rm and}\quad
\La\left(\varphi_n\right)-V\varphi_n=-\lambda_n\,\varphi_n
$$
Choosing $\mu(dx)=\upsilon(dx)/\upsilon(h_0^2)$ and $h_n(x)=\sqrt{\upsilon(h_0^2)}\varphi_n(x)$ we obtain~the formulae (\ref{Psi-spectral-intro-1}) and (\ref{Psi-spectral-intro}) stated in the introduction. The Feynman-Kac propagator version of Corollary~\ref{Poincare-cor} is described below.
\begin{cor}\label{Ka-L2-cor}
For any time horizon $t\geq 0$  and any $f\in \LL_2(\upsilon)$ we have the estimates \begin{equation}\label{expo-decays-Ka}
\left\Vert e^{\lambda_0t}\Ka_{t}(f)-\frac{h_0}{\eta_{\infty}(h_0)}~\eta_{\infty}(f) \right\Vert_{2,\upsilon}\leq e^{-\lambda_1^ht}~\left(\upsilon(f^2)-\frac{\upsilon(h_0)^2}{\upsilon(h_0^2)}~\eta_{\infty}(f)^2
\right)^{1/2}.
\end{equation}
\end{cor}

\medskip

The proof of the above theorem and corollary are given in section~\ref{sec-spectral-proofs}. If, in addition we have $R>0$ and the matrix $A$  is Hurwitz, then we have $\upsilon(1)<\infty$. In this situation, applying (\ref{expo-decays-Ka}) to the unit function $f=1$ yields
\begin{equation}\label{expo-decays-Ka-1}
\left\Vert e^{\lambda_0t}\Ka_{t}(1)-\frac{h_0}{\eta_{\infty}(h_0)} \right\Vert_{2,\upsilon}\leq e^{-\lambda_1^ht}~\left(\upsilon(1)-\frac{\upsilon(h_0)^2}{\upsilon(h_0^2)}
\right)^{1/2}.
\end{equation}
The above results are clearly unsatisfactory when $A$ is not Hurwitz. In this situation, it is preferable to use the non-asymptotic expansions presented in Theorem~\ref{theo-lamda0-K}.

\section{Some illustrations}\label{sec-illustrations}

\subsection{One dimensional models}\label{sec-one-d}
When $r=1$ the reversible condition is trivially met and we have
$$
P_{\infty}=\frac{A+\sqrt{A^2+RS}}{S},\qquad Q_{\infty}=\frac{A+\sqrt{A^2+RS}}{R}\quad \mbox{\rm and}\quad P^h_{\infty}=\frac{R}{2}~\frac{1}{\sqrt{A^2+RS}}.
$$
In this situation, the ground state $h_0$ discussed in (\ref{def-h}) and (\ref{intro-h0}) is given by  
$$
h_0(x):=\exp{\left(-~\frac{A+\sqrt{A^2+RS}}{2R}~x^2\right)}
\quad \mbox{\rm and}\quad 
\lambda_0=\frac{1}{2}\left(A+\sqrt{A^2+RS}\right).
$$
In addition, the eigenfunctions defined in (\ref{intro-def-varphi-ln}) are given,
for any $n\geq 1$, by
$$
\varphi^h_n(x):=\frac{1}{\sqrt{n!}}~
\HH^1_n\left( \sqrt{\frac{2\sqrt{A^2+RS}}{R}}~ x\right)\quad \mbox{\rm and}\quad
\lambda_n^h=n\,\sqrt{A^2+RS}$$
For instance, we have
$$
\varphi^h_1(x)=
 \sqrt{\frac{2\sqrt{A^2+RS}}{R}}~ x\quad \mbox{\rm and}\quad\varphi^h_2(x):=\frac{1}{\sqrt{2}}~\left( \frac{2(A^2+RS)^{1/2}}{R}~ x^2-1\right).
$$
As expected the first excited state $\varphi_1$ of the Hamiltonian is null at the origin, while
$$
\varphi_2(x)=0\Longleftrightarrow x=\pm \sqrt{\frac{R}{2\sqrt{A^2+RS}}}.
$$ 

\subsection{Mehler's formula}\label{Mehler-sec}
We now assume that, $r = 1$, $A=0$ and $R=S^{-1}$. In this situation, we have
$$
Q_{\infty}=S,\qquad
2P^h_{\infty}=P_{\infty}=S^{-1}\quad\mbox{\rm and}\quad
\Lambda^h=-I.
$$
Thus, we may choose $\Za=I$.
In this context, we  readily check that
$$
\upsilon(dx)=dx,\qquad   \lambda=\frac{r}{2}\quad \mbox{\rm and}\quad 
\lambda^{\!V}_n=\sum_{1\leq i\leq r}\left(n_i+\frac{1}{2}~\right)=\vert n\vert+\frac{r}{2}.
$$
The eigenstates are defined by the formulae
 $$
 \varphi_n(x)=\left(\frac{1}{\pi^r}~\mbox{det}\left(S\right)\right)^{1/4}~h_0(x)\varphi^h_n(x)
 $$
 with the re-scaled Hermite polynomials
$$
\varphi^h_n(x):=\frac{1}{\sqrt{n!}}~
\HH^r_n\left(\sqrt{2}~S^{1/2}x\right),
$$
 and the ground state  
$$
h_0(x)=\exp{\left(-\frac{1}{2}\,x^{\prime} S x\right)}.
$$
We also have that
$$
\widehat{X}^h_t(x)=e^{-t}~x
\quad \mbox{\rm and}
\quad
P^h_t=\frac{1}{2}~\left(1-e^{-2t}\right)~S^{-1}~
\Longrightarrow
\mbox{\rm det}(P^h_t)=\frac{\left(1-e^{-2t}\right)^r}{2^r}\frac{1}{\mbox{\rm det}(S)}.
$$
The spectral decomposition reduces to the formula
\begin{align*}
\Ka^h_t(x,dy)
&=\frac{1}{\pi^{r/2}}~\frac{\sqrt{\mbox{\rm det}(S)}}{\left(1-e^{-2t}\right)^{r/2}}~\exp{\left(-\frac{1}{(1-e^{-2t})}~\left(y-e^{-t}~x\right)^{\prime}S\left(y-e^{-t}~x\right)\right)}~dy\\
\\
&=\sum_{n\in\NN^r}~e^{-\vert n\vert t}~\varphi^h_n(x)\varphi^h_n(y)~\eta^h_{\infty}(dy)
\end{align*}
with the collection of functions
$$
\varphi^h_n(x):=\frac{1}{\sqrt{n!}}~
\HH^r_n\left(\sqrt{2}~S^{1/2}x\right)
\quad \mbox{\rm and}\quad
\eta^h_{\infty}=\Na(0,P^h_{\infty})=\Na\left(0,\frac{1}{2}~S^{-1}\right).
$$
Thus, we have
$$
\begin{array}{l}
\displaystyle
\frac{1}{\left(1-e^{-2t}\right)^{r/2}}~\exp{\left(-\frac{1}{(1-e^{-2t})}~\left(y-e^{-t}~x\right)^{\prime}S\left(y-e^{-t}~x\right)\right)}\exp{\left(y^{\prime} Sy\right)}\\
\\
\displaystyle\qquad =\sum_{n\in\NN^r}~e^{-\vert n\vert t}~\frac{1}{\sqrt{n!}}~
\HH^r_n\left(\sqrt{2}~S^{1/2}x\right)\frac{1}{\sqrt{n!}}~
\HH^r_n\left(\sqrt{2}~S^{1/2}y\right).
\end{array}
$$
Replacing $(S^{1/2}x,S^{1/2}y)$ by $(x,y)$ and $e^{-t}$ by $\rho$, we recover Mehler's formula
$$
\begin{array}{l}
\displaystyle
\frac{1}{\left(1-\rho^2\right)^{r/2}}~\exp{\left(-\frac{\rho^2}{1-\rho^2}~(\Vert x\Vert^2+\Vert y\Vert^2)+\frac{2\rho}{1-\rho^2}~ x^{\prime}y\right)}\\
\\
\displaystyle \qquad=\sum_{n\in\NN^r}~\rho^{\vert n\vert}~\frac{1}{\sqrt{n!}}~
\HH^r_n\left(\sqrt{2}~x\right)\frac{1}{\sqrt{n!}}~
\HH^r_n\left(\sqrt{2}~y\right).
\end{array}
$$

The Feynman-Kac propagator takes the form
$$
\Ka_t(x,dy)=\exp{\left(-\frac{rt}{2} \right)}~\exp{\left(-\frac{1}{2}\,x^{\prime} S x\right)}
~K_t^h(x,dy)~\exp{\left(\frac{1}{2}\,y^{\prime} S y\right)}.
$$
On the other hand, we have
$$
\begin{array}{l}
\displaystyle
-\frac{1}{2}\,x^{\prime} S x+\frac{1}{2}\,y^{\prime} S y-\frac{1}{1-e^{-2t}}~\left(y-e^{-t}~x\right)^{\prime}S\left(y-e^{-t}~x\right)\\
\\
\displaystyle \qquad=- y^{\prime}Sy\left(\frac{1}{1-e^{-2t}}-\frac{1}{2}\right)-x^{\prime} S x\left(\frac{1}{2}+\frac{e^{-2t}}{1-e^{-2t}}\right)+2x^{\prime} Sy~\frac{e^{-t}}{1-e^{-2t}}\\
\\
\displaystyle \qquad=-\frac{1}{2}~\left(x^{\prime}Sx+y^{\prime}Sy\right)~ \frac{1+e^{-2t}}{1-e^{-2t}}+x^{\prime} Sy~\frac{2e^{-t}}{1-e^{-2t}}.
\end{array}
$$
Thus, we recover the Mehler's formulation of the Feynman-Kac propagator
$$
\begin{array}{l}
\displaystyle
\Ka_t(x,dy) =\frac{\sqrt{\mbox{\rm det}(S)}}{(2\pi)^{r/2}}~\left(\frac{1}{\sinh(t)}\right)^{r/2}~\exp{\left(-\frac{\coth{(t)}}{2}~\left(x^{\prime}Sx+y^{\prime}Sy\right) +\frac{x^{\prime} Sy}{\sinh(t)}\right)}~dy.
\end{array}$$

\subsection{Quantum harmonic oscillator}\label{multi-dimn-osc-sec}
For diagonal matrices $R$ and $S$ we can choose $\Za=Id$. When $A=0$, the measure $\upsilon$ coincides with the Lebesgue measure $\upsilon(dx)=dx$. 
The $r$-dimensional quantum harmonic oscillator is associated with  diagonal matrices $R$ and $S$  with diagonal entries
 $$
 \begin{array}{l}
\displaystyle
 S_{i,i}=\kappa_i=m\omega_i^2\quad \mbox{\rm and}\quad R_{i,i}=\frac{\hbar^2}{m}\\
 \\
\displaystyle \Longrightarrow \sqrt{\lambda_i(RS)}=\sqrt{(RS)_{i,i}}=
\hbar\omega_i \quad \mbox{\rm and}\quad \left((SR^{-1})_{i,i}\right)^{1/2}=\frac{m\omega_i}{\hbar},
\end{array}
$$
for $1\leq i\leq r$, where $m$ stands for the mass of a particle, $\hbar$ is Planck's constant, and $\omega_i=\sqrt{\kappa_i/m}$ stands for the angular frequencies for some non-negative parameters $\kappa_i$. 
In this situation, have
$$
\lambda^{\!V}_n=\sum_{1\leq i\leq r}\lambda^{\!V}_{n_i}\quad \mbox{\rm and}\quad
 h_n(x)= \prod_{1\leq i\leq r} h_{n_i}(x_i)
$$
with the energy
$$
\lambda^{\!V}_{n_i}:=\left(n_i+\frac{1}{2}~\right)~\hbar\omega_i.
$$
In addition, the eigenfunctions are given by 
$$
\varphi_{n_i}(x_i):=\frac{1}{\sqrt{n_i!}}\left(\frac{m\omega_i}{\pi\hbar}
\right)^{1/4}
~\exp{\left(-\frac{1}{2}\,\frac{m\omega_i}{\hbar}\,
 x_i^2\right)}
~\HH_{n_i}\left(\sqrt{\frac{2m\omega_i}{\hbar}
}
\,x_i\right).
 $$
 The isotropic harmonic oscillator corresponds to the case $\kappa_i=\kappa\Longrightarrow\omega_i=\omega$.
In this case, the energy is given by
$$
\lambda^{\!V}_{n}:=\left(\vert n\vert +\frac{r}{2}~\right)~\hbar\omega.
$$

\subsection{Time varying models}
Theorem~\ref{theo-eta-eta-c}  can be extended to time-varying models associated with time dependent matrices $(A_t,B_t,S_t)$ and $R_t:=B_tB^{\prime}_t$. 
In this situation, the Hamiltonian $\Ha$ in (\ref{def-Schrod-intro}) is a time varying 
operator
$$
\Ha=-\La_t+V_t\quad \mbox{\rm with}\quad V_t(x):=\frac{1}{2}\, x^{\prime}S_tx,
$$
with the second order differential kinetic energy operator
\begin{equation}\label{def-La-intro-t}
\La_t(f)(x):=(A_tx)^{\prime}~\nabla f(x)+\frac{1}{2}~\tr\left(R_t\nabla^2f(x)\right).
\end{equation}
In the above display, $\nabla f$ stands for the gradient column vector with entries $\partial_{x_i} f$, and $\nabla^2 f$ stands for the Hessian matrix with entries $\partial_{x_i,x_j} f$. In the context of multidimensional harmonic oscillators discussed in section~\ref{multi-dimn-osc-sec} we can choose $A=0$ and  diagonal matrices $(R_t,S_t)$  with diagonal entries
$$
(S_t)_{i,i}=\kappa_i=m(t)\omega_i(t)^2\quad \mbox{\rm and}\quad (R_t)_{i,i}=\frac{\hbar^2}{m(t)},
$$
with mass $m(t)$ and angular frequency $\omega_i(t)$.
Replacing $h_0$ by the time varying function
$$
h_t(x)=\exp{\left(-\frac{1}{2}~x^{\prime}Q_{t}\,x\right)}
\quad \mbox{\rm with}\quad
-\partial_tQ_{t}:=
A_t^{\prime}Q_{t}+Q_tA_t- Q_{t}R_tQ_{t} +S_t,
$$
it follows that
$$
h^{-1}_t\left(\partial_t+\La_t\right)(h_t)(x)-V_t(x)=-\lambda_t:=-\tr\left(R_tQ_{t}\right).
$$
Let $\Ca([0,t],\RR^r)$ be the Banach space of all paths from $[0,t]$ to $\RR^r$ equipped with the uniform norm. For any time horizon $t\geq 0$, any measurable function $F_t$ on $\Ca([0,t],\RR^r)$ we find the exponential change of probability formula 
$$
\begin{array}{l}
\displaystyle
\EE\left(F((X_s)_{s\in [0,t]})~\exp{\left(-\int_0^tV_s(X_s)~ds\right)}\right)\\
\\
\displaystyle\qquad=\exp{\left(-\int_0^t\lambda_s ds\right)}~\eta_0(h_0)~\EE\left(F\left((X^{h}_{s})_{s\in [0,t]}\right)~h^{-1}_t(X^h_{t})\right)
\end{array}
$$
with the time varying $h$-process
\begin{equation}\label{def-X-h-t}
dX^h_{t}=\left(A_t-R_tQ_{t}\right)X^h_{t}~dt+B_tdW_t.
\end{equation}

In this context, the free evolution diffusion $X_t$, the $h$-process $X^h_t$ as well as $\widehat{X}_t$ and the Riccati matrix flow $P_t$  are defined as in (\ref{lin-Gaussian-diffusion-filtering}) (\ref{def-X-h-d}) and (\ref{ricc-intro}) by replacing the matrices 
$Q_{\infty}$ by $Q_t$,  and $(A,B)$ and $(R,S)$ by $(A_t,B_t)$ and $(R_t,S_t)$, respectively. For a detailed discussion on time inhomogeneous Riccati equations and related exponential semigroups we refer the reader to~\cite{BishopDelMoralMatricRicc} and references therein. We also mention that~\cite{BishopDelMoralMatricRicc} discusses time-varying controllability conditions that ensures that the stability of the time inhomogeneous version of the evolution equations (\ref{ricc-intro}).

\section{Riccati matrix differential equations}\label{section-Riccati}
\subsection{Gramians  fixed point formulae}\label{gramian-section}
This section is mainly taken from~\cite{bp-21}. 
The difference between the positive and negative fixed points $(P^-_{\infty},P_{\infty})$ of the Riccati equation (\ref{intro-ricc}) is given by the  formula
\begin{equation}
P_{\infty}-P^{-}_{\infty}=\Delta_{\infty}^{-1} \label{def-SS-OO}
\end{equation}
with the Gramian matrices defined via
\begin{equation}\label{gamian-1-ref}
\Delta_t:=
\int_0^{t}\, e^{s(A-P_{\infty}S)^{\prime}}
S~e^{s(A-P_{\infty}S)}\,ds\longrightarrow_{t\rightarrow\infty}\Delta_{\infty}:=\int_0^{\infty}\, e^{s(A-P_{\infty}S)^{\prime}}
S~e^{s(A-P_{\infty}S)}\,ds \in \Sa^+_r.
\end{equation}
Consider now the linear matrix functional
\begin{equation}\label{inter-f-0}
\FF_t~:~P\in \Sa_r^0\mapsto
\FF_t(P):=\left[(\Delta_t^{-1}-\Delta^{-1}_{\infty})+(P-P_{\infty}^-)\right]~\Delta_t~\in \Ga l_r
\end{equation}
Rearranging and using \eqref{gamian-1-ref} implies that
\begin{equation}\label{inter-f-1}
\FF_t(P)=I+(P-P_{\infty})\,\Delta_t \quad \text{ and } \quad \FF_t(P_{\infty})=I.
\end{equation}
Recall that $Q_{\infty}$ is defined as $P_{\infty}$ by replacing $(A,R,S)$ by $(A^{\prime},S,R)$.
In the same vein, $\Delta_t^h$ is defined as $\Delta_t$ by replacing $(A,R,S)$ and $P_{\infty}$ by $(A^{\prime},S,R)$ and $Q_{\infty}$. Thus, by symmetry arguments and (\ref{formula-example-Pinfty-0-Q}), we also have
\begin{equation}\label{Delta-P-Q}
(\Delta_{\infty}^h)^{-1}=P^{-1}_{\infty}-(P^-_{\infty})^{-1}=Q_{\infty}-Q_{\infty}^-\quad \mbox{\rm and}\quad \left(\phi^h_{t}(0),P^h_{\infty}\right)=\left(\Delta_t^h,\Delta_{\infty}^h\right),
\end{equation}
 with the Gramian matrices
$$
\Delta_t^h:=
\int_0^t  e^{(A-RQ_{\infty})s}~R~e^{(A-RQ_{\infty})^{\prime}s}~ds\longrightarrow_{t\rightarrow\infty}\Delta_{\infty}^h:=\int_0^{\infty}\,  e^{(A-RQ_{\infty})s}~R~e^{(A-RQ_{\infty})^{\prime}s}~ds\in \Sa^+_r.
$$

The following lemma proves the second equality on the right hand side of \eqref{intro-h0}.
\begin{lem}\label{lem-QP-2}
For any $(A,R,S)$ satisfying the rank condition (\ref{def-contr-obs}), we have
\begin{equation} \label{key-QP-2}
\tr(SP_{\infty})=2~\tr(A)-\tr(SP^{-}_{\infty})=\tr(Q_{\infty}R).
\end{equation}
\end{lem}
\proof
The Gramian $\Delta_{\infty}$ satisfies the Sylvester equations given by
$$
(A-P_{\infty}S)^{\prime}\Delta_{\infty}+\Delta_{\infty}(A-P_{\infty}S)+S=0=\Delta_{\infty}^{-1}(A-P_{\infty}S)^{\prime}+(A-P_{\infty}S)\Delta_{\infty}^{-1} +\Delta_{\infty}^{-1} S\,\Delta_{\infty}^{-1}.
$$
It then follows that
$$
\begin{array}{l}
\displaystyle
\tr\left(\Delta_{\infty}\left(\Delta_{\infty}^{-1}(A-P_{\infty}S)^{\prime}+(A-P_{\infty}S)\Delta_{\infty}^{-1}+\Delta_{\infty}^{-1}S\,\Delta_{\infty}^{-1}
\right)\right)=0\\
\\
\Longrightarrow 2\tr(SP_{\infty})=2\tr(A)+\tr(S\,\Delta_{\infty}^{-1})=2\tr(A)+\tr(SP_{\infty})-\tr(SP^{-}_{\infty}),
\end{array}
$$
where we have used \eqref{def-SS-OO} to obtain the final equality.
From this we obtain
\begin{equation} \label{key-QP}
\tr(SP^{-}_{\infty})=2\tr(A)-\tr(SP_{\infty}).
\end{equation}
In the same vein, we have
\begin{eqnarray*}
 \tr\left(P^{-1}_{\infty}\mbox{\rm Ricc}(P_{\infty})\right)=0&\Longrightarrow&
\tr(SP_{\infty})= 2\tr(A)+\tr(P^{-1}_{\infty}R)\nonumber\\
 \tr\left((P^-_{\infty})^{-1}\mbox{\rm Ricc}(P^-_{\infty})\right)=0&\Longrightarrow&
\tr(SP^-_{\infty})= 2\tr(A)+\tr((P^-_{\infty})^{-1}R).
\end{eqnarray*}
 Combining the last assertion with (\ref{key-QP}) we conclude that
 $$
\tr(P_{\infty}S)=- \tr((P^-_{\infty})^{-1}R)=\tr(Q_{\infty}R),
 $$
as required.
 \cqfd

\subsection{A Floquet-type representation} \label{ref-floquet-sec}
For any $P\in \Sa^0_r$ and $\delta>0$ set
 $$
\chi(P):= \Vert P_{\infty}^-\Vert^{-1}
\left[\Vert P_{\infty}-P_{\infty}^-\Vert+\Vert P-P_{\infty}\Vert\right]\quad\mbox{and}\quad
\chi_{\delta}:=\left[{\lambda_{ min}\left(\Delta_{\delta}\right)
\lambda_{ min}\left(-P_{\infty}^-\right)}\right]^{-1}.
$$
The next theorem provides an explicit description of $\Ea_t(P)$ in terms of the matrices $(A,S,P_{\infty})$.
 
\begin{theo}[Floquet-type representation~\cite{bp-21}]\label{theo-floquet}
For any time horizon $t\geq 0$ and any $P\in\Sa_r^0$ we have Riccati exponential  semigroup formula
\begin{equation}\label{floquet-type-ric}
\Ea_{t}(P)= e^{t (A-P_{\infty}S) }~\FF_t(P)^{-1} = \Ea_t(P_\infty)\FF_t(P)^{-1},
\end{equation}
where $\FF_t(P)$ was defined in (\ref{inter-f-0}).
For any $t\geq \delta>0$  we have the uniform estimates
\begin{eqnarray}
{\Vert \FF_t(P)^{-1}\Vert} \,\leq \chi_{\delta}&\mbox{and}&
\Vert \Ea_t(P)\Vert \,\leq~\chi_{\delta}\,\Vert \Ea_t(P_{\infty})\Vert.
\label{estimate-E-noQ}
\end{eqnarray}
\textcolor{black}{In addition, for any $t\geq 0$ we have the exponential estimates \eqref{exp-stab-phi} as well as the bounds}
\begin{eqnarray}
{\Vert \FF_t(P)^{-1}\Vert}\,\leq \chi(P)&\mbox{and}&~~ \quad
\Vert \Ea_t(P)\Vert \,\leq~\chi(P)\,\Vert \Ea_t(P_{\infty})\Vert.
\label{estimate-EQ}
\end{eqnarray}
\end{theo}
Using the decomposition
\begin{equation}
\mbox{\rm Ricc}(Q_1)-\mbox{\rm Ricc}(Q_2)
\,=\,(A-Q_1S)(Q_1-Q_2)+(Q_1-Q_2)(A-Q_2S)^{\prime},
\label{polarization-formulae}
\end{equation}
for $Q_1, Q_2 \in \Sa_r^0$, \textcolor{black}{applying (\ref{floquet-type-ric})} we have the closed form Lipschitz type matrix formula
\begin{equation}\label{eqn-cor-decomp}
\phi_t(Q_1)-\phi_t(Q_2)=\Ea_t(P_{\infty})~\FF_t(Q_1)^{-1}(Q_1-Q_2)~\left(\Ea_t(P_{\infty})~\FF_t(Q_2)^{-1}\right)^{\prime}.
\end{equation}
Applying (\ref{eqn-cor-decomp}) with $Q_2=P_{\infty}$ and using (\ref{inter-f-1}), we recover the Bernstein-Prach-Tekinalp formula~\cite{prach2015infinite,prach-thesis} given by
\begin{equation}\label{pbt-ref}
\phi_t(P)=P_{\infty}+ \Ea_t(P_{\infty})~\FF_t(P)^{-1}(P-P_{\infty})~\Ea_t(P_{\infty})^{\prime}.
\end{equation}

\subsection{Lipschitz inequalities}\label{ref-sect-expo-cv}
Combining Theorem~\ref{theo-floquet}, \eqref{exp-stab-phi} and \eqref{eqn-cor-decomp} we easily obtain the following result.
\begin{theo}\label{cor-phi-Lip}
For any time horizon 
$t\geq \delta>0$ and  any $Q_1,Q_2\in\Sa_r^0$ we have the Lipschitz estimate
$$
\Vert \phi_t(Q_1)-\phi_t(Q_2)\Vert\leq  (\alpha\chi_{\delta})^2\,e^{-2\beta t}~\Vert Q_1-Q_2\Vert
$$
with the parameters $(\alpha,\beta, \chi_{\delta})$ defined in (\ref{exp-stab-phi}) and Theorem~\ref{theo-floquet}.
In addition, for any $t\geq 0$ we have the local Lipschitz estimate
$$
\Vert \phi_t(Q_1)-\phi_t(Q_2)\Vert\leq  \alpha^2\chi(P_1)\chi(P_2)~e^{-2\beta t}~\Vert Q_1-Q_2\Vert,
$$
with the parameters $\chi(Q_i)$ defined in Theorem~\ref{theo-floquet}.
\end{theo}
Noting that
\begin{eqnarray*}
\Ea_t(Q_1)-\Ea_t(Q_2)&=&\Ea_t(P_{\infty})~\FF_t(Q_1)^{-1}~\left[\FF_t(Q_2)-\FF_t(Q_1)\right]~\FF_t(Q_2)^{-1}\\
&=&\Ea_t(P_{\infty})~\FF_t(Q_1)^{-1}~(Q_2-Q_1)~\Delta_t~\FF_t(Q_2)^{-1},
\end{eqnarray*}
where $\Delta_t$ was defined in \eqref{gamian-1-ref}, we also obtain the following corollary.
\begin{cor}
For any time horizon 
$t\geq \delta>0$ and  any $Q_1,Q_2\in\Sa_r^0$ we have the Lipschitz estimate
$$
\Vert \Ea_t(Q_1)-\Ea_t(Q_2)\Vert\leq  \alpha\,\chi_{\delta}^2\,\Vert\Delta_{\infty}\Vert\,e^{-\beta t}~\Vert Q_1-Q_2\Vert
$$
with the parameters $(\alpha,\beta, \chi_{\delta})$ defined in (\ref{exp-stab-phi}) and Theorem~\ref{theo-floquet}.
In addition, for any $t\geq 0$ we have local Lipschitz estimate
$$
\Vert \Ea_t(Q_1)-\Ea_t(Q_2)\Vert\leq \alpha\Vert\Delta_{\infty}\Vert\,\chi(Q_1)\chi(Q_2)\,e^{-\beta t}~\Vert Q_1-Q_2\Vert
$$
with the parameter $\chi(Q_i)$ defined in Theorem~\ref{theo-floquet}.
\end{cor}

The first coordinate of the evolution semigroup  (\ref{coupled-equations}) can be written as
$$
\widehat{X}_t(x,P_0)=\Ea_t(P_0) x
$$
Using the decomposition
$$
\widehat{X}_t(x_1,Q_1)-\widehat{X}_t(x_2,Q_2)=(\Ea_t(Q_1)-\Ea_t(Q_2)) x_1+\Ea_t(Q_2) (x_1-x_2),
$$
we readily check the following theorem.
\begin{theo}\label{cor-X-hat-lip}
For any time horizon 
$t\geq \delta>0$ and  any $Q_1,Q_2\in\Sa_r^0$ we have the  estimate
$$
\Vert \widehat{X}_t(x_1,Q_1)-\widehat{X}_t(x_2,Q_2)\Vert\leq  \alpha\chi_{\delta}\,e^{-\beta t}~ \left(\chi_{\delta}\Vert\Delta_{\infty}\Vert\,\,\Vert x_1\Vert~\Vert Q_1-Q_2\Vert+\Vert x_1-x_2\Vert\right)
$$
with the parameters $(\alpha,\beta, \chi_{\delta})$ defined in (\ref{exp-stab-phi}) and Theorem~\ref{theo-floquet}.
In addition, for any $t\geq 0$ we have the estimate
$$
\Vert\widehat{X}_t(x_1,Q_1)-\widehat{X}_t(x_2,Q_2)\Vert\leq \ \alpha\,\chi(Q_2)\,e^{-\beta t}~ \left(\chi(Q_1)\Vert\Delta_{\infty}\Vert\,\,\Vert x_1\Vert~\Vert Q_1-Q_2\Vert+\Vert x_1-x_2\Vert\right)
$$
with the parameter $\chi(Q_i)$ defined in Theorem~\ref{theo-floquet}.
\end{theo}

\section{Feynman-Kac propagators}
\subsection{Unnormalised semigroups}

Whenever the initial state $X_0$ is distributed according to some probability measure $\eta_0$ on $\RR^d$, by Fubini's theorem and (\ref{def-Ka}) we have
\begin{eqnarray*}
\eta_0(\psi_t)&:=&\int
\eta_0(dx)\psi_t(x)=\int \eta_0(dx)\int\Ka_t(x,dy)\psi_0(y)\\
&=&\int\left(\int \eta_0(dx)\Ka_t(x,dy)\right)\psi_0(y)=\left(\eta_0\Ka_t\right)(\psi_0).
\end{eqnarray*}
This yields the formula
$$
\eta_0(\psi_t)=\gamma_{t}(\psi_0)\quad \mbox{\rm with the measure}\quad 
\gamma_{t}:=\eta_0\Ka_t.
$$
Equivalently, the measure $\gamma_t$ is defined by the unnormalised Feynman-Kac path integral
$$
\gamma_t(\psi_0)=\EE\left(\psi_0(X_t)~\exp{\left(-\int_0^t V(X_u)~ds\right)}\right)\quad \mbox{\rm and}\quad \gamma_t(1)=\EE\left(\exp{\left(-\int_0^t V(X_u)~ds\right)}\right).
$$
In the above display, $1$ stands for the unit function.
Observe that the evolution semigroup of $\gamma_t$ is linear and given by the formulae
\begin{equation}\label{sg-prop}
 \gamma_{t+s}=\eta_0\Ka_{s+t}=\eta_0\left(\Ka_{s}\Ka_t\right)=\left(\eta_0\Ka_s\right)\Ka_t=\gamma_s\Ka_t.
\end{equation}	
 Finally observe that 
$$
\partial_t\gamma_t(\psi_0)=\EE\left(\La(\psi_0)(X_t)~\exp{\left(-\int_0^t V(X_u)~ds\right)}\right)-\EE\left(\psi_0(X_t)~V(X_t)~\exp{\left(-\int_0^t V(X_u)~ds\right)}\right)
$$
This yields the evolution equation
\begin{equation}\label{evol-gamma}
\partial_t\gamma_t(\psi_0)=\gamma_t(\La(\psi_0))-\gamma_t(\psi_0V)=-\gamma_t(\Ha(\psi_0)),
\end{equation}
where $\Ha$ was defined in \eqref{def-H-intro}.

\subsection{Normalised semigroups}\label{sec-norm-measures}
We shall denote by $\eta_t$ the normalised probability measures
\begin{equation}\label{gamma-1}
\eta_t(f):=\gamma_t(f)/\gamma_t(1)\Longrightarrow\gamma_t(f)=\eta_t(f)~\exp{\left(-\int_0^t\eta_s(V)~ds\right)}.
\end{equation}	
We check this claim using the formula
$$
-\partial_t\log{\gamma_t(1)}=\frac{1}{\gamma_t(1)}~\EE\left(V(X_t)~\exp{\left(-\int_0^t V(X_u)~ds\right)}\right)=\eta_t(V).
$$
In bra-ket notation, the semigroup property (\ref{sg-prop}) and the probability measure $\eta_t$ can be written in the form
$$
\eta_0=\mu_{\varphi}\Longrightarrow \langle \varphi\vert e^{-(s+t)\Ha}\vert \psi\rangle=
 \langle \varphi\vert e^{-s\Ha}e^{-t\Ha}\vert \psi\rangle\quad \mbox{\rm and}\quad
 \eta_t(\psi_0)=\frac{ \langle \varphi~\vert e^{-t\Ha}~\vert \psi_0\rangle}{ \langle \varphi~\vert e^{-t\Ha}~\vert 1\rangle}.
$$
Observe that for any $s\leq t$ we have the correspondence principle
\begin{equation}
\eta_t(f)=\frac{\gamma_t(f)}{\gamma_t(1)}=\frac{\gamma_s\Ka_{t-s}(f)}{\gamma_s\Ka_{t-s}(1)}=\frac{\eta_s\Ka_{t-s}(f)}{\eta_s\Ka_{t-s}(1)}=:\Phi_{t-s}(\eta_s)(f).
\label{def-Phi-semi}
\end{equation}
The semigroup $\Phi_{t-s}(\eta_s)=\eta_t$ of the normalised measures described above is a nonlinear mapping from the set of probability measures on $\RR^r$ into itself.

\subsection{Normalised Feynman-Kac propagators}\label{sec-norm-propagators}

There two different ways to normalise the integral Feynman-Kac operators $\Ka_t$. The first one is based on (\ref{gamma-1}), which implies that
\begin{equation}\label{unnorm-key}
\eta_t(f)=\gamma_t(f)\exp{\left(\int_0^t\eta_s(V)~ds\right)}.
\end{equation}
This yields the formula
$$
\eta_t(f)=
\EE\left(f(X_t)~\exp{\left(-\int_0^t V_{\eta_s}(X_s)~ds\right)}\right)\quad \mbox{\rm with}\quad V_{\eta_s}(x)=V(x)-\eta_s(V).
$$
This shows that the normalised measures  $\eta_t$ are defined as $\gamma_t$ by replacing $V$ by the time varying centered potential function $V_{\eta_t}$. 
It is therefore natural to consider the normalised propagator defined below.
\begin{defi}
For any initial distribution $\eta_0$ and for any $t\geq 0$ we denote by $\Ka^{\eta_0}_{t}$ the integral operator
$$
\Ka^{\eta_0}_{t}(f)(x)
:=\EE\left(f(X_t)~\exp{\left(-\int_0^t V_{\eta_s}(X_s)~ds\right)}~|~X_0=x\right).
$$
\end{defi}
Using \eqref{def-Phi-semi}, it is straightforward to see that
\begin{eqnarray}
\Ka^{\eta_0}_{t}(f)(x)
&=&\exp{\left(\int_0^t\Phi_s(\eta_0)(V)~ds\right)}\times
\Ka_t(f)(x)=\frac{\Ka_t(f)(x)}{\eta_0\Ka_t(1)}. 
\label{eq-Ka-1}
\end{eqnarray}
From this we deduce that
$$
\eta_0\Ka^{\eta_0}_{t}(f)=\frac{1}{\gamma_t(1)}~\eta_0\Ka_t(f)=\frac{\gamma_t(f)}{\gamma_t(1)}\Longrightarrow \eta_0\Ka^{\eta_0}_{t}=\eta_t.
$$
\begin{prop}
For any $s,t\geq 0$ we have the evolution semigroup properties
$$
\Ka^{\eta_0}_{s+t}=\Ka^{\eta_0}_{s}\,\Ka^{\eta_s}_{t},\qquad
\eta_{s+t}=\eta_s\Ka^{\eta_s}_{t}\quad \mbox{and}\quad 
\Ka^{\eta_{\infty}}_{t}(f)=e^{\lambda_0 t}~\Ka_t(f)
$$
\end{prop}
\proof
To prove the two semigroup properties, note that
\begin{eqnarray*}
\Ka^{\eta_0}_{t}(f)(x)
&=&\exp{\left(\int_0^s\eta_s(V)~ds\right)}~\exp{\left(\int_s^t\Phi_{u-s}(\eta_s)(V)~du\right)}
~\Ka_{t+s}(f)(x)\\
&=&\exp{\left(\int_0^{t-s}\Phi_{u}(\eta_s)(V)~du\right)}~\Ka^{\eta_0}_{s}(\Ka_t(f))(x)=\Ka^{\eta_0}_{s}\left(\Ka^{\Phi_{s}(\eta_0)}_{t-s}(f)\right)(x).
\end{eqnarray*} 
The second then follows from $\eta_0\Ka^{\eta_0}_{t}=\eta_t$. 

To check the right hand side, note that
$$
\eta_{\infty}(\Ka_t(1))=\exp{\left(-\int_0^t\Phi_s(\eta_{\infty})(V)~ds\right)}=e^{-\eta_{\infty}(V)t}=e^{-\lambda_0t}.
$$
The result then follows from \eqref{eq-Ka-1}. \cqfd

Arguing as in (\ref{evol-gamma}) with $V$ replaced by $V_{\eta_t}$, we also find the evolution equation
\begin{equation}
\overline{\psi}_t:=\Ka^{\eta_0}_{t}(\psi_0)\Longrightarrow
-\partial_t\overline{\psi}_t=\Ha_{\eta_t}(\overline{\psi}_t)\quad \mbox{\rm and}\quad
-\partial_t\eta_t(f)=\eta_t\left(\Ha_{\eta_t}(f)\right),
\label{evol-eta}
\end{equation}
with the normalised time varying Hamiltonian
$$
\Ha_{\eta_t}=-\La+(V-\eta_t(V))\Longleftrightarrow (\ref{evol-eta-cov}).
$$

\medskip

A second strategy to normalise the Feynman-Kac propagator is to divide by its total mass.  
\begin{defi}\label{def-overline-K}
We associate with $\Ka_t$ the normalised the Markov integral operator
 $ \overline{\Ka}_t$  defined by the ratio formula
$$
 \overline{\Ka}_t(f)(x):={\Ka_t(f)(x)}/{\Ka_t(1)(x)}=\int~f(y)~\Phi_t(\delta_x)(dy).
$$
\end{defi}
By (\ref{sg-prop}) and (\ref{gamma-1}) the normalising constant is given by
$$
\Ka_t(1)(x)=\delta_x \Ka_t(1)=\exp{\left(-\int_0^t\Phi_s(\delta_x)(V)~ds\right)}.
$$
Using Theorem~\ref{theo-intro-entrop-wasserstein} (see also the estimates presented in section~\ref{ref-sect-expo-cv}) it follows that
$$
\begin{array}{l}
\delta_x=\Na(x,0)\\
\\
\displaystyle\Longrightarrow
 \overline{\Ka}_t(x,dy)=\Phi_t(\delta_x)(dy)=\Na\left(\widehat{X}_t(x,0),\phi_t(0)\right)(dy)\longrightarrow_{t\rightarrow\infty}
 \eta_{\infty}=\Na(0,P_{\infty}).
 \end{array}
$$
This yields the formula
$$
2\Phi_t(\delta_x)(V)=\widehat{X}_t(x,0)^{\prime}S\widehat{X}_t(x,0)+\tr(S\phi_t(0))=x^{\prime}\Ea_t(0)^{\prime}S\Ea_t(0)x+\tr(S\phi_t(0)),
$$
from which we conclude that
$$
\Ka_t(1)(x)=\exp{\left(-\frac{1}{2}x^{\prime}\Delta_t(0)x-\frac{1}{2}\int_0^t\tr(S\phi_s(0))ds\right)},
$$
with
$$
 \Delta_t(0):=\int_0^t\Ea_s(0)^{\prime}S\,\Ea_s(0)~ds\longrightarrow_{t\rightarrow\infty}
  \Delta_{\infty}(0):=\int_0^{\infty}\Ea_s(0)^{\prime}S\,\Ea_s(0)~ds.
$$
In contrast with linear semigroups we have the nonlinear transport formula
$$
\Phi_t(\eta_0)=\BB_{\Ka_t(1)}(\eta_0)\overline{\Ka}_t\not= \eta_0\overline{\Ka}_t,
$$
and thus, the following theorem.
\begin{theo}\label{theo-non-Gauss}
For any probability measure $\eta$ on $\RR^r$ and for any $t\geq 0$ we have the Boltzmann-Gibbs formula
$$
\Phi_t(\eta)=\BB_{\varpi_t}(\eta)\overline{\Ka}_t
$$
with the energy function
$$
\log{\varpi_t(x)}=-\frac{1}{2}\,x^{\prime}\Delta_t(0) \,x\longrightarrow_{t\rightarrow \infty}
\log{\varpi_{\infty}(x)}=-\frac{1}{2}\,x^{\prime}\Delta_{\infty}(0) \,x.
$$
In addition, for any $\delta>0$ there exists some constant $c_{\delta}$ such that
for any pair of probability measures, $\eta$ and $\mu$, on $\RR^r$ and for $t\geq \delta>0$, we have
\begin{equation}\label{ref-non-Gauss}
 \Vert\Phi_t(\eta)-\Phi_t(\mu)\Vert_{\tiny tv}\leq \frac{c_{\delta}}{\eta(\varpi_{\infty})\wedge \mu(\varpi_{\infty})}~e^{-\beta t},
\end{equation}
where the parameter $\beta$ was introduced in (\ref{exp-stab-phi}).
\end{theo}

\proof
{\color{black}The first assertion comes from the fact that
$$
\Ka_t(1)(x)=\varpi_t(x)~\exp{\left(-\frac{1}{2}\int_0^t\tr(S\phi_s(0))ds\right)}\Longrightarrow
\BB_{\Ka_t(1)}=\BB_{\varpi_t}.
$$}
To prove \eqref{ref-non-Gauss}, observe that
$$
\Phi_t(\eta)(f)-\Phi_t(\mu)(f)=\int~\BB_{\varpi_t}(\eta)(dx)\BB_{\varpi_t}(\mu)(dy)
\left(\overline{\Ka}_t(f)(x)-\overline{\Ka}_t(f)(y)\right).
$$
From Definition \ref{def-overline-K} and Theorem~\ref{theo-intro-entrop-wasserstein}, for any function $f$ such that $\Vert f\Vert\leq 1$, we have
$$
\vert \overline{\Ka}_t(f)(x)-\overline{\Ka}_t(f)(y)\vert\leq \Vert\Phi_t(\delta_x)-\Phi_t(\delta_y)\Vert_{\tiny tv} \leq c_{\delta}~e^{-\beta t}~\Vert x-y\Vert.
$$
This implies that
$$
 \Vert\Phi_t(\eta)-\Phi_t(\mu)\Vert_{\tiny tv} \leq 2c_{\delta}~e^{-\beta t}~
\left( \int~\BB_{\varpi_t}(\eta)(dx)\Vert x\Vert\vee  \int~\BB_{\varpi_t}(\mu)(dx)\Vert x\Vert\right).
$$
On the other hand, we have
$$
\eta(\varpi_t)\geq \eta(\varpi_{\infty})>0\quad \mbox{\rm with}\quad \log{\varpi_{\infty}(x)}=-\frac{1}{2}\,x^{\prime}\Delta_{\infty}(0) \,x,
$$
as well as the uniform estimate
$$
\int~\eta(dx) \varpi_t(x)~\Vert x\Vert\leq \int~\eta(dx) \varpi_{\delta}(x)~\Vert x\Vert\leq \sup_{x\in\RR^r}\vert \varpi_{\delta}(x)~\Vert x\Vert\vert =: c_\delta'<\infty, \quad t \ge \delta > 0.
$$
Thus
$$
\int~\BB_{\varpi_t}(\eta)(dx)\Vert x\Vert = \frac{1}{\eta(\varpi_t)}\int\eta(dx) \varpi_t(x)~\Vert x\Vert \le \frac{c_\delta'}{\eta(\varpi_{\infty})},
$$
which completes the proof of the theorem. 
\cqfd

\medskip

Observe that
$$
\eta=\Na(x,P)\Longrightarrow \left\{
\begin{array}{rcl}
\BB_{\varpi_t}(\eta)&=&\Na\left(\left(P^{-1}+\Delta_t(0)\right)^{-1}P^{-1}x,\left(P^{-1}+\Delta_t(0)\right)^{-1}\right)\\
&&\\
\eta \overline{\Ka}_t&=&\Na\left(\Ea_t(0)x,\phi_t(0)+\Ea_t(0)P\Ea_t(0)^{\prime}\right).
\end{array}\right.$$
Combining  these two formula we readily check the following proposition.
\begin{prop}\label{alternative-description-xw-P}
For any $t\geq 0$, we have the Gaussian preserving property
$$
\eta_0=\Na(x,P)\Longrightarrow \eta_t=\Na\left(\widehat{X}_t(x,P),\phi_t(P)\right)
$$
with the mean vector and covariance matrices given by
$$
\left\{
\begin{array}{rcl}
\widehat{X}_t(x,P)&=&\Ea_t(0)\left(P^{-1}+\Delta_t(0)\right)^{-1}P^{-1}x\\
&&\\
\phi_t(P)&=&\phi_t(0)+\Ea_t(0)\left(P^{-1}+\Delta_t(0)\right)^{-1}\Ea_t(0)^{\prime}
\end{array}\right.
$$
\end{prop}

\subsection{Some stability properties}
\subsubsection{Relative entropy}\label{entropy-sec}

For any initial conditions of the form given in (\ref{spec-initial-Gauss}),
the Boltzmann-Kullback Leibler relative entropy of $\Phi_t(\eta_0)$ with respect to $\Phi_t(\mu_0)$ has a closed form (see for instance the article~\cite{zhang} and formula (A.23) in~\cite{rasmussen}) given by the formula
\begin{eqnarray}
\displaystyle
 &&\mbox{\rm Ent}\left(\Phi_t(\eta_0)~|~\Phi_t(\mu_0)\right)\notag\\
 \notag \\
\displaystyle &&\qquad=\frac{1}{2}\left(
 \mbox{\tr}\left(\phi_t(P)\phi_t(Q)^{-1}-I\right)+\log{\mbox{det}\left(\phi_t(Q)\phi_t(P)^{-1}\right)}\right.
 \notag\\
 \notag\\
\displaystyle&&\hskip2cm \left.+\left(\widehat{X}_t(x,P)-\widehat{X}_t(y,Q)\right)^{\prime}\phi_t(Q)^{-1}\left(\widehat{X}_t(x,P)-\widehat{X}_t(y,Q)\right)\right).
\label{eq-ent-BKL}
 \end{eqnarray}
 When $P=0=Q$ the above formula reduces to
 $$
 \mbox{\rm Ent}\left(\Phi_t(\delta_x)~|~\Phi_t(\delta_y)\right)=\frac{1}{2}(x-y)^{\prime} \Ea_t(0)^{\prime}\phi_t(0)^{-1}\Ea_t(0)(x-y).
 $$
 To estimate determinants of matrices close to the identity we use  the following lemma.
\begin{lem}[\cite{tugaut-dm-1}]\label{lem-tech-2}
For any $(r\times r)$-matrix $A$  we have
$$
\Vert A \Vert_F< \frac{1}{2}\Longrightarrow
\left\vert\log{\mbox{\rm det}\left(I-A\right)}\right\vert 
\leq \frac{3}{2}~\Vert A \Vert_F~.
$$
\end{lem}
 
\bigskip 

For any $n\geq 1$ and $\delta>0$ we  set
$$
t_{n,\delta}:= \delta\vee\frac{1}{2\beta}
\log{\left(2n \alpha_{\delta}\right)}~\quad \mbox{\rm with}\quad \alpha_{\delta}:=(\alpha\chi_{\delta})^2\left(r\lambda_{\rm max}(\Pi_{-,\delta}^{-1})~\tr(\Pi_{-,\delta}^{-1})\right)^{1/2}.
$$
with the positive matrix $\Pi_{-,\delta}$ and the parameters $(\alpha,\beta, \chi_{\delta})$ defined in (\ref{exp-stab-phi}) and Theorem~\ref{theo-floquet}. 
\begin{theo}\label{theo-time-delta-n}
For any initial conditions
$$
\eta_0=\Na(x,P)\quad \mbox{and}\quad\mu_0=\Na(y,Q) \quad\mbox{s.t.}\quad \Vert P-Q\Vert\leq n,
$$
and for any $ t\geq t_{\delta,n}$ with $n\geq 1$ and $\delta>0$, we have the exponential decay estimate
$$
\begin{array}{l}
\displaystyle
 \mbox{\rm Ent}\left(\Phi_t(\eta_0)~|~\Phi_t(\mu_0)\right)\\
\\
\displaystyle\leq \left(\frac{5}{4}~\alpha_{\delta}~\Vert P-Q\Vert+\lambda_{\rm max}(\Pi_{-,\delta}^{-1})~ (\alpha\chi_{\delta})^2\,~ \left((\chi_{\delta}\Vert\Delta_{\infty}\Vert)^2\,\,\Vert x\Vert^2~\Vert P-Q\Vert^2+\Vert x-y\Vert^2
\right)\right)~e^{-2\beta t}.
\end{array}$$
\end{theo}
 \proof
We start by controlling the $\log{\rm det}$ term on the right hand side of \eqref{eq-ent-BKL} using Lemma \ref{lem-tech-2}.
First note that, by (\ref{exp-stab-phi}), for any $ t\geq \delta >0$ we have
 $$
 0<\Pi_{+,\delta}^{-1}\leq \phi_t(Q)^{-1}\leq \Pi_{-,\delta}^{-1}\Longrightarrow 
 \tr(\phi_t(Q)^{-2})\leq r\lambda_{\rm max}(\Pi_{-,\delta}^{-1})~\tr(\Pi_{-,\delta}^{-1}).
 $$
 Combining this with the Lipschitz estimates stated in Theorem~\ref{cor-phi-Lip} we check that
 $$
\Vert I-\phi_t(Q)^{-1}\phi_t(P)\Vert_F\leq \left(r\lambda_{\rm max}(\Pi_{-,\delta}^{-1})~\tr(\Pi_{-,\delta}^{-1})\right)^{1/2}~(\alpha\chi_{\delta})^2\,e^{-2\beta t}~\Vert P-Q\Vert.
 $$
{\color{black}For any $t\geq t_{n,\delta}$ we have
$$
\begin{array}{l}
\displaystyle t\geq \delta\quad \mbox{\rm and}\quad e^{-2\beta t}\leq\frac{1}{
2n \alpha_{\delta}}~\quad \mbox{\rm with}\quad \alpha_{\delta}:=(\alpha\chi_{\delta})^2\left(r\lambda_{\rm max}(\Pi_{-,\delta}^{-1})~\tr(\Pi_{-,\delta}^{-1})\right)^{1/2}\\
\\
\displaystyle\Longrightarrow
\Vert I-\phi_t(Q)^{-1}\phi_t(P)\Vert_F\leq \frac{1}{2n}~\Vert P-Q\Vert.
\end{array}$$}
Applying Lemma~\ref{lem-tech-2} to $A=\phi_t(Q)^{-1}(\phi_{t}(Q)-\phi_t(P))$, for any 
$t\geq t_{n,\delta}$ and
 $\Vert P-Q\Vert\leq n$ we have the estimate
  \begin{eqnarray*}
\left\vert\log{\mbox{det}\left(\phi_t(P)\phi_t(Q)^{-1}\right)}\right\vert
&\leq &\frac{3}{2}~\Vert \phi_t(P)-\phi_t(Q) \Vert_F~\Vert \phi_t(Q)^{-1}\Vert_F\\
&\leq &\frac{3}{2}~\left(r\lambda_{\rm max}(\Pi_{-,\delta}^{-1})~\tr(\Pi_{-,\delta}^{-1})\right)^{1/2}~(\alpha\chi_{\delta})^2\,e^{-2\beta t}~\Vert P-Q\Vert.
\end{eqnarray*}
Similarly, for any $t\geq \delta$ we have
  \begin{eqnarray*}
\left\vert \mbox{\rm tr}\left(I-\phi_t(Q)^{-1}\phi_t(P)\right)\right\vert&\leq&  \left\Vert \phi_t(Q)^{-1}\right\Vert_F~\left\Vert\phi_t(P)-\phi_t(Q)\right\Vert_F\\
&\leq &\left(r\lambda_{\rm max}(\Pi_{-,\delta}^{-1})~\tr(\Pi_{-,\delta}^{-1})\right)^{1/2}~\left\Vert\phi_t(P)-\phi_t(Q)\right\Vert\\
&\leq &\left(r\lambda_{\rm max}(\Pi_{-,\delta}^{-1})~\tr(\Pi_{-,\delta}^{-1})\right)^{1/2}
 (\alpha\chi_{\delta})^2\,e^{-2\beta t}~\Vert P-Q\Vert.
\end{eqnarray*}
Finally, we notice that
$$
\left\vert
\left(\widehat{X}_t(x,P)-\widehat{X}_t(y,Q)\right)^{\prime}\phi_t(Q)^{-1}\left(\widehat{X}_t(x,P)-\widehat{X}_t(y,Q)\right)
\right\vert\leq \lambda_{\rm max}(\Pi_{-,\delta}^{-1})~ \Vert (\widehat{X}_t(x,P)-\widehat{X}_t(y,Q)\Vert^2.
$$
Applying Theorem~\ref{cor-X-hat-lip} for any $t\geq \delta$ we check that
$$
\begin{array}{l}
\displaystyle
\left\vert
\left(\widehat{X}_t(x,P)-\widehat{X}_t(y,Q)\right)^{\prime}\phi_t(Q)^{-1}\left(\widehat{X}_t(x,P)-\widehat{X}_t(y,Q)\right)
\right\vert\\
\\
\displaystyle\leq \lambda_{\rm max}(\Pi_{-,\delta}^{-1})~ (\alpha\chi_{\delta})^2\,e^{-2\beta t}~ \left(\chi_{\delta}\Vert\Delta_{\infty}\Vert\,\,\Vert x\Vert~\Vert P-Q\Vert+\Vert x-y\Vert
\right)^2,
\end{array}
$$
which concludes the proof.
\cqfd
\medskip

Applying the above theorem to $P=0$ and $(y,Q)=(0,P_{\infty})$ we check the following corollary
\begin{cor}
For any $ t\geq t_{\delta,n}$ with $n=\lfloor \Vert P_{\infty}\Vert\rfloor $ and $\delta>0$, we have the exponential decay estimate
$$
\begin{array}{l}
\displaystyle
 \mbox{\rm Ent}\left(\Phi_t(\delta_x)~|~\eta_{\infty}\right)\leq  \left(\frac{5}{4}~\alpha_{\delta}~\Vert P_{\infty}\Vert+\lambda_{\rm max}(\Pi_{-,\delta}^{-1})~ (\alpha\chi_{\delta})^2\,~ \left(1+(\chi_{\delta}\Vert\Delta_{\infty}\Vert~\Vert P_{\infty}\Vert)^2\right)~\Vert x\Vert^2\right)~e^{-2\beta t}.
 \end{array}
 $$
In addition, for any $t\geq \delta>0$ and any $x,y\in\RR^r$ we have
$$
 \mbox{\rm Ent}\left(\Phi_t(\delta_x)~|~\Phi_t(\delta_y)\right)\leq \frac{1}{2}\lambda_{\rm max}(\Pi_{-,\delta}^{-1})~ (\alpha\chi_{\delta})^2\,e^{-2\beta t}~ \Vert x-y\Vert^2.
 $$
\end{cor}

\subsubsection{Wasserstein distances}\label{wasserstein-sec}
For any initial conditions of the form (\ref{spec-initial-Gauss}) we have
$$
\WW_2\left(\Phi_t(\mu_0),\Phi_t(\eta_0)\right)^2\leq \Vert \widehat{X}_t(x,P)-\widehat{X}_t(y,Q)\Vert^2+\Vert\phi_t(P)^{1/2}-\phi_t(Q)^{1/2}\Vert_F^2
$$
For any $P,Q\in\Sa_{r}^+$ we also have the Ando-Hemmen inequality
\begin{equation}\label{square-root-key-estimate}
\Vert P^{1/2}- Q^{1/2}\Vert \leq \left[\lambda^{1/2}_{min}(P)+\lambda^{1/2}_{min}(Q)\right]^{-1}~\Vert P-Q\Vert
\end{equation}
for any unitary invariant matrix norm $\Vert . \Vert$ (including the spectral and the Frobenius norms). See for instance Theorem 6.2 on page 135 in~\cite{higham}, as well as Proposition 3.2  in~\cite{hemmen}. For a more thorough discussion on the geometric properties of positive semidefinite matrices and square roots we refer to~\cite{hiriart}.

Using (\ref{exp-stab-phi}) and theorem~\ref{cor-phi-Lip} for any $t\geq \delta>0$ we check that
$$
\Vert\phi_t(P)^{1/2}-\phi_t(Q)^{1/2}\Vert_F\leq  \sqrt{r}\left[2\lambda^{1/2}_{min}(\Pi_{-,\delta})\right]^{-1}~ (\alpha\chi_{\delta})^2\,e^{-2\beta t}~\Vert P-Q\Vert
$$
Using theorem~\ref{cor-X-hat-lip} we obtain the following theorem
\begin{theo}\label{ref-W-intro}
For any $t\geq \delta>0$ and any initial conditions
$$
\eta_0=\Na(x,P)\quad \mbox{and}\quad\mu_0=\Na(y,Q) 
$$
we have the exponential decay estimate
$$
\begin{array}{l}
\displaystyle
\WW_2\left(\Phi_t(\mu_0),\Phi_t(\eta_0)\right)\\
\\
\displaystyle\leq \alpha\chi_{\delta}\,e^{-\beta t}\left(~\Vert x-y\Vert+\chi_{\delta}\left( \Vert\Delta_{\infty}\Vert\,\,\Vert x\Vert+
\alpha  \sqrt{r}\left[2\lambda^{1/2}_{min}(\Pi_{-,\delta})\right]^{-1}~\,e^{-\beta t}\right)~\Vert P-Q\Vert\right)
\end{array}
$$

\end{theo}

\subsubsection{Proof of Theorem~\ref{theo-lamda0-K}}\label{theo-lamda0-K-proof}
Observe that
$$
\Ka_t(h_0)(x)=e^{-\lambda_0 t}~h_0(x)\Longrightarrow
e^{\lambda_0t}\Ka_t(1)(x)=\frac{h_0(x)}{\overline{\Ka}_t(h_0)(x)}=\frac{h_0(x)}{\eta_{\infty}(h_0)}~k_t(x)
$$
with
$$
k_t(x):={\eta_{\infty}(h_0)}/{\overline{\Ka}_t(h_0)(x)}={\eta_{\infty}(h_0)}/{\Phi_t(\delta_x)(h_0)}
$$
Observe that
$$
k_t(x)
=\left(\frac{\mbox{\rm det}(I+\phi_t(0)Q_{\infty})}{\mbox{det}\left(I+P_{\infty}Q_{\infty}\right)}\right)^{1/2}
\exp{\left(\frac{1}{2}~\widehat{X}_t(x,0)^{\prime}(I+Q_{\infty}\phi_t(0))^{-1}\widehat{X}_t(x,0)\right)}
$$
This shows that
$$
k_t(x)\longrightarrow_{\vert x\vert\rightarrow+\infty}+\infty\quad \mbox{\rm and}\quad 
k_t(x)\longrightarrow_{t\rightarrow\infty}1
$$
This shows that one cannot expect uniform upper bound with respect to the state variable. 
The estimates (\ref{estim-final}) are now a direct consequences of the Lipschitz estimates presented
in section~\ref{ref-sect-expo-cv}, which we now detail.

Applying (\ref{pbt-ref}) to $P=0$ we check the formula
$$
(I+\phi_t(0)Q_{\infty})(I+P_{\infty}Q_{\infty})^{-1}=I- \Ea_t(P_{\infty})~\FF_t(0)^{-1}P_{\infty}~\Ea_t(P_{\infty})^{\prime}(Q_{\infty}^{-1}+P_{\infty})^{-1}
$$
On the other hand, by (\ref{estimate-E-noQ}) for any $t\geq \delta>0$  we have
$$
\Vert \Ea_t(P_{\infty})~\FF_t(0)^{-1}P_{\infty}~\Ea_t(P_{\infty})^{\prime}(Q_{\infty}^{-1}+P_{\infty})^{-1}\Vert\leq \chi_{\delta}^3~\Vert \Ea_t(P_{\infty})\Vert^2~\Vert P_{\infty}(Q_{\infty}^{-1}+P_{\infty})^{-1}~\Vert
$$
Using (\ref{exp-stab-phi}) this implies that
$$
\Vert \Ea_t(P_{\infty})~\FF_t(0)^{-1}P_{\infty}~\Ea_t(P_{\infty})^{\prime}(Q_{\infty}^{-1}+P_{\infty})^{-1}\Vert_F\leq \sqrt{r}~ \chi_{\delta}^3~\alpha^2\,e^{-2\beta t}~\Vert P_{\infty}Q_{\infty}(I+P_{\infty}Q_{\infty})^{-1}~\Vert
$$
Applying lemma~\ref{lem-tech-2} for any 
$$
t>t_{\delta}:=\frac{1}{2\beta}\log{\left(
2\sqrt{r}~ \chi_{\delta}^3~\alpha^2\,~\Vert P_{\infty}Q_{\infty}(I+P_{\infty}Q_{\infty})^{-1}~\Vert\right)}
$$
we have
$$
\left\vert\log{\left(\frac{\mbox{\rm det}(I+\phi_t(0)Q_{\infty})}{\mbox{det}\left(I+P_{\infty}Q_{\infty}\right)}\right)^{1/2}}\right\vert 
\leq \frac{3}{4}~\sqrt{r}~ \chi_{\delta}^3~\alpha^2\,e^{-2\beta t}~\Vert P_{\infty}Q_{\infty}(I+P_{\infty}Q_{\infty})^{-1}~\Vert
$$
On the other hand, we have
$$
\widehat{X}_t(x,0)^{\prime}(I+Q_{\infty}\phi_t(0))^{-1}\widehat{X}_t(x,0)\leq \Vert \phi_t(0)^{-1}\Vert\,\Vert(\phi_t(0)^{-1}+Q_{\infty})^{-1}\Vert~\Vert \widehat{X}_t(x,0)\Vert^2
$$
By (\ref{exp-stab-phi}) for any $t\geq \delta>0$ we have
$$
 \phi_t(0)^{-1}\leq \Pi^{-1}_{-,\delta}\quad \mbox{\rm and}\quad
 (\phi_t(0)+Q_{\infty})^{-1}\leq (\Pi_{-,\delta}+Q_{\infty})^{-1}
$$
combining the above with the estimates (\ref{exp-stab-phi}) and (\ref{estimate-E-noQ}) we check that
$$
\widehat{X}_t(x,0)^{\prime}(I+Q_{\infty}\phi_t(0))^{-1}\widehat{X}_t(x,0)\leq \Vert \Pi^{-1}_{-,\delta}\Vert\,\Vert(\Pi_{-,\delta}+Q_{\infty})^{-1}\Vert~(\alpha\chi_{\delta})^2\,~e^{-2\beta t}~\Vert x\Vert^2
$$

To finish the proof of the theorem, note that for any time horizon $t\geq 0$ and any $f\in \LL_1(\eta_{\infty})$ we have the decomposition
\begin{eqnarray*}
e^{\lambda_0t}~\Ka_t(f)(x)
&=&\frac{h_0(x)}{\eta_{\infty}(h_0)}~\left(\eta_{\infty}(f)+\widetilde{K}_t(f)(x)\right)~k_t(x)
\end{eqnarray*}
with the integral operator
$$
\widetilde{K}_t(f)(x)=\overline{K}_t(f)(x)-\eta_{\infty}(f)=\Phi_t(\delta_x)(f)-\eta_{\infty}(f)
$$
This ends the proof of theorem~\ref{theo-lamda0-K}.\cqfd

\section{Path integral formulations}
\subsection{Particle absorption models}\label{particle-absorption-sec}

For any time horizon $t\geq 0$, any measurable function $F_t$ on $\Ca([0,t],\RR^r)$ and any starting point $x\in\RR^r$, we have the Feynman-Kac path-integral formula
$$
\EE\left(F_t\left((X_s)_{s\in [0,t]}\right)~\exp{\left(-\int_0^t V(X_u)~ds\right)}~|~X^c_0=x\right)=\EE(F_t\left((X^c_s)_{s\in [0,t]}\right)~1_{\tau^c\geq t}~|~X_0=x).
$$ 
We recover the Feynman-Kac propagator  formula (\ref{def-Ka}) by choosing
functions that only depend on the terminal time.

\begin{theo}\label{theo-c-h0}
For all $x \in \RR^r$, we have
\begin{equation}\label{La-h0}
h_0^{-1}\La(h_0)(x)=V(x)-\lambda_0,
\end{equation}
where $h_0$ was defined in \eqref{def-h}.
In addition, for any time horizon $t\geq 0$, any measurable function $F_t$ on $\Ca([0,t],\RR^r)$ and any starting point $x\in\RR^r$, we have the Feynman-Kac path-integral formula
\begin{equation}\label{F-h0}
\begin{array}{l}
\displaystyle
\EE(F_t\left((X^c_s)_{s\in [0,t]}\right)~1_{\tau^c\geq t}~|~X_0=x)\\
\\
\displaystyle\quad=\exp{\left(-\lambda_0 t\right)}~h_0(x)~\EE\left(F_t\left((X^{h}_{s})_{s\in [0,t]}\right)~h_0^{-1}(X^h_{t})~|~X^h_0=x\right),
\end{array}
\end{equation}
where $X^{h}_{t}$ stands for the diffusion with generator defined by
$$
\La^{h}(f)=\La(f)+h_0^{-1}\Gamma_{\La}(h_0,f),
$$
with the carr\'e-du-champ operator
$$
\Gamma_{\La}(h_0,f)(x)
:=\left(B\nabla h_0(x)\right)^{\prime}\left(B\nabla f(x)\right)=-h_0 (x)~(RQ_{\infty}x)^{\prime}~\nabla f(x).
$$

\end{theo}

\proof
From \eqref{def-h}, we have the gradient formula
$$
\begin{array}{l}
\displaystyle\nabla \log{h_0(x)}=-Q_{\infty}x\quad \mbox{\rm and}\quad \nabla^2\log{h_0(x)}=-Q_{\infty}\\
\\
\displaystyle\Longrightarrow h_0^{-1}(x)\nabla h_0(x)=-Q_{\infty}x\quad \mbox{\rm and}\quad
 h_0^{-1}(x)\nabla^2h_0(x)=-Q_{\infty}+(Q_{\infty}x)\left(Q_{\infty}x\right)^{\prime}.
\end{array}$$
This implies that
\begin{eqnarray*}
h^{-1}_0\La(h_0)(x)&=&h_0^{-1}(x)\left((Ax)^{\prime}\nabla h_0(x)+\frac{1}{2}\tr\left(R\nabla^2h_0\right)\right)\\
&=&-x^{\prime} A^{\prime}Q_{\infty}x+\frac{1}{2}\tr\left(R(Q_{\infty}x)\left(Q_{\infty}x\right)^{\prime}\right)-\frac{1},{2}\tr\left(RQ_{\infty}\right)
\end{eqnarray*}
from which it follows that
\begin{eqnarray*}
h_0^{-1}\La(h_0)(x)-V(x)&=&-\frac{1}{2}x^{\prime} (A^{\prime}Q_{\infty}+Q_{\infty}A)x+\frac{1}{2}~x^{\prime}(Q_{\infty}RQ_{\infty})x-\frac{1}{2}x^{\prime}Sx-\frac{1}{2}\tr\left(RQ_{\infty}\right)\\
&=& -\lambda_0,
\end{eqnarray*}
where the last equality follows from the fact that $A^{\prime}Q_{\infty}+Q_{\infty}A-Q_{\infty}RQ_{\infty}+S=0$, as in \eqref{def-h}.

This yields the exponential change of probability formula 
$$
\begin{array}{l}
\displaystyle
\exp{\left(\lambda_0 t\right)}~\EE\left(F((X_s)_{s\in [0,t]})~\exp{\left(-\int_0^tV(X_s)~ds\right)}\right)\\
\\
\displaystyle=\EE\left(F((X_s)_{s\in [0,t]})~\frac{h_0(X_0)}{h_0(X_t)}~\left(\frac{h_0(X_t)}{h_0(X_0)}\exp{\left(-\int_0^t(h_0^{-1}\La h_0)(X_s)~ds\right)}\right)\right)\\
\\
\displaystyle=~\eta_0(h_0)~\EE\left(F\left((X^{h}_{s})_{s\in [0,t]}\right)~h_0^{-1}(X^h_{t})\right),
\end{array}
$$
which ends the proof of the theorem.
\cqfd

Thus, combining this with Lemma \ref{lem-QP-2}, Theorem \ref{theo-1-intro} is now proved.
Moreover, Theorem~\ref{theo-eta-eta-c} is a now a direct consequence of (\ref{F-h0}). Indeed, using (\ref{F-h0})
we verify that the law $\eta^h_t$ of the random states $X^h_t$ of the $h$-process defined in (\ref{def-X-h-d}) and the distribution of the non-absorbed particle defined in (\ref{def-d-non-abs}) are connected by the Boltzmann-Gibbs transformation (\ref{BG-h0-eta-intro}).
The Gaussian preserving property of the linear diffusion process $X^h_t$ is immediate.
The formulae given in (\ref{ref-link-h-wX}) are easily checked using the the Boltzmann-Gibbs transformation (\ref{BG-h0-eta-intro}). Moreover, we can check that the pair $(\widehat{X}_t,P_t)$ 
given by (\ref{ref-link-h-wX}) satisfies (\ref{ricc-intro}) using brute force calculations, or 
by checking that the Gaussian distributions with mean and covariance matrices $(\widehat{X}_t,P_t)$ solves the nonlinear equation (\ref{evol-eta-cov}).

\subsection{Backward $h$-processes}\label{section-backward}

For a fixed time horizon $t\geq 0$, we let  $\overline{X}_t$ be a random sample from $\Na(\widehat{X}_t,P_t)$. We also  denote by
$X^h_{t,s}(x)$, with $s\in [0,t]$, be the backward diffusion defined by
$$
dX^h_{t,s}(x)=\left(A X^h_{t,s}(x)+RP_s^{-1}( X^h_{t,s}(x)-\widehat{X}_s)\right)~ds+BdW_s,
$$
starting at $X^h_{t,t}(x)=x$ at time $s=t$. In the above display, $P_s$ stands for the solution of the Riccati matrix differential equation defined in (\ref{ricc-intro}).  We assume that $\overline{X}_t$ and $(W_s)_{s\leq t}$ are independent.

Rewritten in terms of the density $g_s$ of the Gaussian distribution $\Na(\widehat{X}_s,P_s)$, we have
\begin{equation}\label{theo-eq}
 \overline{X}_{t,s}^h:=X^{h}_{t,s}(\overline{X}_t)
\Longrightarrow d \overline{X}_{t,s}^h=\left(A \overline{X}_{t,s}^h-R\,\nabla \log{g_s( \overline{X}_{t,s}^h)}\right)\,ds+
B\,dW_s.
\end{equation}

The following theorem, taken from \cite{bdc-21} links the non-absorbed particle process with the above backward diffision

\begin{theo}[\cite{bdc-21}]\label{theo-backward}
Assume that $X_0\sim \Na(\widehat{X}_0,P_0)$. In this situation, for any $t\geq 0$ we have the backward formulation of the Feynman-Kac path integral
$$
\begin{array}{l}
\displaystyle
\EE\left(F\left((X_s^c)_{s\in [0,t]}\right)~|~\tau^c\geq t\right)=\EE\left(F\left((\overline{X}^h_{t,s})_ {s\in [0,t]}\right)\right).\end{array}
$$

\end{theo}
The random state $\overline{X}_{t,s}^h$ is a Gaussian variable with a mean $\widehat{X}^h_{t,s}$ and covariance matrix $P^h_{t,s}$
satisfying the backward equations
$$
\left\{
\begin{array}{rcl}
\partial_s\widehat{X}^h_{t,s}&=&A \widehat{X}^h_{t,s}+RP_s^{-1}(\widehat{X}^h_{t,s}-\widehat{X}_s)\\
&&\\
\partial_{s}P_{t,s}^h&=&(A+RP^{-1}_s)P^h_{t,s}+P^h_{t,s}(A+RP^{-1}_s)^{\prime}-R
\end{array}\right.$$
with the terminal condition $(\widehat{X}^h_{t,t},P_{t,t}^h)=(\widehat{X}_t,P_t)$, where $(\widehat{X}_s,P_s)$ is the solution to the forward equations described in (\ref{ricc-intro}).

%
%

\subsection{Extensions to nonlinear diffusions}

The $h$-process methodology can be extended to more general generators $\La$ and other choices of the potential function $V$. 
We now assume that $\La$ is the generator of the diffusion equation
\begin{equation}\label{nonlin-diffusion}
dX_t=A(X_t)dt+B(X_t)dW_t
\end{equation}
for some drift function $A(x)$ and some diffusion matrix valued function $B(x)$ with appropriate dimensions. We also assume there exists some ground state $h_0$ associated with some energy $\lambda_0$; that is, we have that
$$
h_0^{-1}\La(h_0)(x)=V(x)-\lambda_0
$$
In this situation, the $h$-process $X^{h}_{t}$ is a diffusion with generator defined by
$$
\La^{h}(f)=\La(f)+h_0^{-1}\Gamma_{\La}(h_0,f)$$
with the carr\'e-du-champ operator
$$
\Gamma_{\La}(h_0,f
)(x):=\left(B(x)\nabla h_0(x)\right)^{\prime}\left(B(x)\nabla f(x)\right)=
\left(R(x)\nabla h_0(x)\right)^{\prime}\nabla f(x),
$$
where we have defined $R(x):=B(x)B(x)^{\prime}$. Equivalently, the $h$-process  is defined by the diffusion 
$$
dX^{h}_{t}=\left(A(X^{h}_{t})+R(X^{h}_{t})\nabla \log{h_0(X^{h}_{t})}\right) dt+B(X^{h}_{t})dW_t
$$

Let  $\overline{X}_t$  a random sample from the Feynman-Kac probability measures $\eta_t$ defined as in (\ref{gamma-1})  for some potential function $V$.  

Whenever it exists,  let $g_s$ be the density of the normalised or unnormalised  Feynman-Kac
measures $\eta_s$ or $\gamma_s$.
In this situation, following the analysis developed in~\cite{bdc-21}, the assertion of Theorem~\ref{theo-backward} remains valid with the backward diffusion
\begin{equation}\label{theo-eq-nl}
  d \overline{X}_{t,s}^h=\left(A(\overline{X}_{t,s}^h)-\mbox{\rm div}_{R}\log{g_s( \overline{X}_{t,s}^h)}\right)\,ds+
B(\overline{X}_{t,s}^h)\,dW_s
\end{equation}
with the terminal condition $\overline{X}_{t,t}^h=\overline{X}_{t}$ and the
$R$-divergence $m$-column vector operator with $j$-th entry given by the formula
\begin{equation*}
\mbox{\rm div}_{R}(f)(x)^j:=~\sum_{1\leq i\leq r}~\partial_{x_i}\left(R_{i,j}(x)~f(x)\right).
\end{equation*}

   \section{McKean-Vlasov interpretations}\label{mcKean-sec}
\subsection{Interacting jump processes}\label{IPS-section}
Let $\overline{X}_t$ be a nonlinear jump diffusion process with generator
$$
\La_{\overline{\eta}_t}(f)(x)=\La(f)(x)+V(x)~\int~(f(y)-f(x))~\overline{\eta}_t(dy)\quad \mbox{\rm where}\quad\overline{\eta}_t := \mbox{\rm Law}(\overline{X}_t).
$$
The process starts at $\overline{X}_0=X_0$.
Between the jumps the process $\overline{X}_t$ evolves as $X_t$. At rate $V(\overline{X}_t)$ the process jumps onto a new location randomly selected according to the distribution $\overline{\eta}_t$. Observe that 
\begin{eqnarray*}
\partial_t\overline{\eta}_t(f)
&=&\overline{\eta}_t\left(\La_{\overline{\eta}_t}(f)\right)=\overline{\eta}_t(\La(f))-\overline{\eta}_t(fV)+\overline{\eta}_t(f)\overline{\eta}_t(V).
\end{eqnarray*}
This shows that $\overline{\eta}_t$ satisfies the same evolution equation as the one satisfied by $\eta_t$ given in (\ref{evol-eta-cov}).
Thus, for any choice of the generator $\La$ and any choice of the  potential function $V$ we have that
$$
\overline{\eta}_t(dx)=\eta_t(dx):=\PP(X^c_t\in dx~|~\tau^c>t).
$$
The mean field particle interpretation of the nonlinear process $\overline{X}_t$ is defined by
 a system of $N$ walkers, $\xi^i_t$, evolving independently as $X_t$ with jump rate $V_t(X_t)$, for $1\leq i\leq N$.
At each jump time, the particle $\xi^{i}_t$ jumps onto a particle uniformly chosen in the pool. 
 The occupation measure of system is given by the empirical measure
 \begin{equation}\label{dmc-1}
\eta^N_t=\frac{1}{N}\sum_{1\leq i\leq N}\delta_{\xi^i_t}\longrightarrow_{N\rightarrow\infty}\eta_t\longrightarrow_{t\rightarrow\infty}\eta_{\infty}.
\end{equation}
Mimicking (\ref{unnorm-key}) we also define the normalising constant approximations
 \begin{equation}\label{dmc-2}
\frac{1}{t}\int_0^t \eta^N_s(V) ds:=-\frac{1}{t}\log{\gamma^N_t(1)}\longrightarrow_{N\rightarrow\infty}-\frac{1}{t}\log{\gamma_t(1)}\longrightarrow_{t\rightarrow\infty}.\lambda_0=\eta_{\infty}(V)
\end{equation}
The interacting particle system discussed above belongs to the class of diffusion Monte Carlo algorithms, see for instance the series of articles~\cite{caffarel,cances-tony,tony,tony-rs,mathias,mathias-2}, as well as~\cite{dm-2000-moran,dm-2000,dmpenev,dm-sch} and the references therein. 
Observe that the $N$ ancestral lines $\zeta_t^i:=(\xi^i_t)_{0\leq s\leq t}$ of length $t$ of the above genetic-type process can also be seen as a 
system of $N$ path-valued particles evolving independently as the historical process $Y_t=(X_{s})_{0\leq s\leq t}$ of $X_t$, with a jump rate $V_t(X_t)$ that only depends on the terminal state $X_t$ of the ancestral line $Y_t$.

\subsection{Interacting diffusions}\label{EnKF-sec}
For any probability measure $\eta$ on $\RR^r$ we let $\Pa_{\eta}$ denote the $\eta$-covariance
\begin{equation}
\eta\mapsto\Pa_{\eta}:=\eta\left([e-\eta(e)][e-\eta(e)]'\right)
\end{equation}
where $e(x):=x$ is the identity function and $\eta(f)$ is a column vector whose $i$-th entry is given by $\eta(f^i)$ for some measurable function $f:\RR^r\rightarrow\RR^r$.

We now consider three different nonlinear McKean-Vlasov-type diffusion process,
\begin{equation}\label{Kalman-Bucy-filter-nonlinear-ref}
\begin{array}{lrcl}
(1)& d\overline{X}_t~&=&~(A-\Pa_{\overline{\eta}_t}S)\,\overline{X}_t~dt+\Pa_{\overline{\eta}_t}~S^{1/2}~
d\Wa_{t}+B\,d\overline{\Wa}_t,\\
&&&\\
(2)& \displaystyle d\overline{X}_t~&=&\displaystyle~\left(A\,\overline{X}_t~-\frac{1}{2}~\Pa_{\overline{\eta}_t}\,S\left(\overline{X}_t+
\overline{\eta}_t(e)\right)\right)dt+B\,d\overline{\Wa}_t,\\
&&&\\
(3)&  \partial_t\overline{X}_t~&=&\displaystyle
A\,\overline{X}_t~-\frac{1}{2}~\Pa_{\overline{\eta}_t}\,S\left(\overline{X}_t+
\overline{\eta}_t(e)\right)+
(R+M_t)\,\Pa_{\overline{\eta}_t}^{-1}\left(\overline{X}_t-\overline{\eta}_t(e)\right),
\end{array}
\end{equation}
for any skew symmetric matrix $M^\prime_t=-M_t$ that may also depend $\overline{\eta}_t$.  In all three cases  $(\Wa_t,\overline{\Wa}_t)$ are independent copies of $W_t$; and $\overline{X}_0$ is an independent copies of $X_0$.  We also assume that $(\Wa_t,\Wa_t,\overline{X}_0)$ are independent. 
 In all three cases in (\ref{Kalman-Bucy-filter-nonlinear-ref}), $\overline{\eta}_t $ stands for the probability distribution of $\overline{X}_t$; that is,  we have that
\begin{equation}\label{def-nl-cov}
\overline{\eta}_t := \mbox{\rm Law}(\overline{X}_t).
\end{equation}
 Observe that, in all three cases the stochastic processes discussed above depend in some nonlinear fashion on the law of the diffusion process itself. 
 \begin{theo}
 In all the three cases presented in (\ref{Kalman-Bucy-filter-nonlinear-ref}), for any $t\geq 0$ we have the Gaussian preserving property
$$
	\eta_0=\Na(\widehat{X}_0,P_0)=\overline{\eta}_0\Longrightarrow
\overline{\eta}_t=\Na(\widehat{X}_t,P_t)=\eta_t.
$$
\end{theo}
\proof
Let $\overline{X}_t$ be the process defined as in $(1)$ by replacing $\Pa_{\overline{\eta}_t}$ by $P_t$.
In this case, we have
$$
d\left(\overline{X}_t-\EE(\overline{X}_t)\right)~=~(A-P_tS)\,\left(\overline{X}_t-\EE(\overline{X}_t)\right)~dt+P_t~S^{1/2}~d\Wa_{t}+B\,d\overline{\Wa}_t.\
$$
Applying Ito's formula and taking expectations we obtain
$$
\partial_t\Pa_{\overline{\eta}_t}=(A-P_tS)\Pa_{\overline{\eta}_t}+\Pa_{\overline{\eta}_t}(A-P_tS)^{\prime}+P_tSP_t+R.
$$
This yields the linear system
$$
\partial_t\left(\Pa_{\overline{\eta}_t}-P_t\right)=(A-P_tS)(\Pa_{\overline{\eta}_t}-P_t)+\Pa_{\overline{\eta}_t}(A-P_tS)^{\prime}(\Pa_{\overline{\eta}_t}-P_t)\Longrightarrow
P_t=\Pa_{\overline{\eta}_t}.
$$
We conclude that $\overline{X}_t$ is a linear diffusion with mean $\widehat{X}_t$ and covariance matrix $P_t$. The proof for the other two cases follows the same lines of arguments, thus we leave the details to the reader.
\cqfd

The mean-field particle interpretation
of the first nonlinear diffusion process in \eqref{Kalman-Bucy-filter-nonlinear-ref} is given by the Mckean-Vlasov type interacting diffusion process
\begin{equation}\label{fv1-3}
\begin{array}{rcl}
d\xi^i_t&=&(A-P^N_tS)\,\overline{X}_t~dt+P^N_t~S^{1/2}~
d\Wa_{t}^i+B\,d\overline{\Wa}^i_t, \qquad i=1,\ldots,N,
\end{array}
\end{equation}
where  $(\Wa_{t}^i,\overline{\Wa}^i_t,\xi^i_0)_{1\leq i\leq N}$  are
 $N$ independent copies of $(\Wa_{t},\overline{\Wa}_t\overline{X}_0)$. 
In the above display,  the $P^N_t$ are the rescaled empirical covariance matrices given by the formulae
\begin{equation}\label{fv1-3-2}
P^N_t:=\left(1-\frac{1}{N}\right)^{-1}~\Pa_{\eta^{N}_t}=\frac{1}{N-1}\sum_{1\leq i\leq N}\left(\xi^i_t-m^N_t\right)\left(\xi^i_t-m^N_t\right)^{\prime}\,,
\end{equation}
with the empirical measures
$$
\eta^{N}_t:=\frac{1}{N}\sum_{1\leq i\leq N}\delta_{\xi^i_t}
\quad\mbox{\rm and the sample mean}\quad
 m^N_t:=\frac{1}{N}\sum_{1\leq i\leq N}\xi_t^i\,.
$$

Note that (\ref{fv1-3}) is a set of $N$ stochastic differential equations coupled via the empirical
covariance matrix $P_t^N$. The mean-field particle interpretation
of the second and third nonlinear diffusion processes in \eqref{Kalman-Bucy-filter-nonlinear-ref} are defined as above by replacing $\Pa_{\eta_t}$ by the sample covariance matrices $P^N_t$.
The  quasi-invariant measure $\eta_{\infty}$ and the parameter $\lambda_0$ are computed using the limiting  formulae (\ref{dmc-1}) and (\ref{dmc-2}).

The interacting diffusions discussed above belong to the class of Ensemble Kalman filters, see for instance the pioneering article by Evensen~\cite{evensen03}, the series of articles~\cite{ap-2016,BishopDelMoralMatricRicc,Bishop/DelMoral/Bruno2012},  as well as~\cite{tugaut-dm-1,tugaut-dm-2}  and the references therein.

In contrast with the interacting jump process discussed in section~\ref{IPS-section} none of the nonlinear diffusions discussed in (\ref{Kalman-Bucy-filter-nonlinear-ref}) can be extended to  more general generators $\La$ and other choices of the  potential function $V$.

We end this section with an application of the seminal feedback particle filter methodology recently developed by Mehta and Meyn and their co-authors~\cite{prashant-1,prashant-2,prashant-3,prashant-4,prashant-5} to Feynman-Kac models. Consider the diffusion 
$$
d\overline{X}_t=\left(A(\overline{X}_t)+U_t(\overline{X}_t)\right)dt+B(\overline{X}_t)dW_t,
$$
where $U_t(x)$ is the solution of the Poisson equation
$$
\sum_{1\leq i\leq r}\frac{1}{g_t(x)}\partial_{x_i}\left(U^i_t(x)~g_t(x)\right)=(V(x)-\eta_t(V)), \qquad t \ge 0.
$$
In the above display $g_t(x)$ stands for the density of the distribution $\overline{\eta}_t$ of the random state $\overline{X}_t$.
The generator $\La_{\overline{\eta}_t}$ of the above time varying diffusion satisfies the equation
$$
\overline{\eta}_t(\La_{\overline{\eta}_t}(f))=\overline{\eta}_t(\La(f))+\sum_{1\leq i\leq r}\int~U^i_t(x) \partial_{x_i}f(x)~g_t(x)~dx.
$$
Integrating by part the last term we obtain the formula
$$
\overline{\eta}_t(\La_{\overline{\eta}_t}(f))=\overline{\eta}_t(L(f))-\int~f(x)~(V(x)-\overline{\eta}_t(V))~\overline{\eta}_t(dx),
$$
from which we conclude that
$$
\overline{\eta}_t(\La_{\overline{\eta}_t}(f))=\overline{\eta}_t(L(f))-\overline{\eta}_t(fV)+\overline{\eta}_t(f)\overline{\eta}_t(V).
$$
This shows that $\overline{\eta}_t=\mbox{\rm Law}(\overline{X}_t)=\eta_t$ coincides with the normalised Feynman-Kac measures. 

For linear-Gaussian models we have $\eta_t=\Na(\widehat{X}_t,P_t)$. Thus,  the Poisson equation resumes to the formula
\begin{eqnarray*}
 \sum_{1\leq i\leq r}\partial_{x_i}U^i_t(x)-(x-\widehat{X}_t)^{\prime}P_t^{-1}U_t(x)&=&\frac{1}{2}(x^{\prime}Sx-\widehat{X}_t^{\prime}S\widehat{X}_t-\tr(SP_t))\\
 &=&\frac{1}{2}~(x-\widehat{X}_t)^{\prime}S(x-\widehat{X}_t).
\end{eqnarray*}
The solution of the above equation is clearly given by
$$
U_t(x)=-\frac{1}{2}~P_tS(x+\widehat{X}_t)\Longrightarrow
 \sum_{1\leq i\leq r}\partial_{x_i}U^i_t(x)=\tr(P_tS).
$$
The resulting diffusion coincides with the second case in (\ref{Kalman-Bucy-filter-nonlinear-ref}).

\section{Spectral decompositions}
The main focus of this section is the proof of Theorem \ref{spectral-theo-intro} and Theorem \ref{theo-Ka-spec}, and the corresponding corollaries. Thus, in what follows, we assume that $R>0$ and $ A^{\prime}=R^{-1}AR$. Before giving the proofs, we first spend some time discussing some properties of the reversible $h$-process introduced in section \ref{rev-models-intro} and the Chebychev-Hermite polynomials introduced in section \ref{spectral-intro-sec}.

{\color{black}\subsection{Reversible $h$-processes}\label{rev-section}
Due to the reversibility conditions, the fixed points of the algebraic Riccati equation (\ref{intro-ricc}) are given by the formulae   
\begin{equation}\label{ref-Weyls-2-intro}
(P_{\infty}^-)^{-1}=-R^{-1}A-R^{-1}(A^2+RS)^{1/2}< 0<
 P_{\infty}^{-1}=-R^{-1}A+R^{-1}(A^2+RS)^{1/2},
\end{equation}
with the square root $(A^2+RS)^{1/2}$ that has all positive eigenvalues. A proof of the above result can be found in~\cite{bp-21}. To check that this square root is well-defined, observe that
\begin{equation}\label{ref-positive-A2}
A^{\prime}=R^{-1}AR\Longrightarrow
A^2+RS=R(A^{\prime}R^{-1} A+S)=R((A^{\prime})^2 +SR)R^{-1},
\end{equation}
which has positive eigenvalues.
We also have the formulae
  \begin{equation}\label{ref-Weyls-2-rev-intro}
  Q_{\infty}=P_{\infty}^{-1}+2 R^{-1}A\quad \mbox{\rm and}\quad A-RQ_{\infty}=-(A^2+RS)^{1/2},
\end{equation}
which yields
$$
\begin{array}{l}
\displaystyle \mbox{\rm Spec}( A-RQ_{\infty})=\left\{-\vert \lambda\vert^{1/2}~:~\lambda\in\mbox{\rm Spec}(A^2+RS) \right\}=\{-\lambda^h_1,\ldots,-\lambda^h_r\}\subset \RR_-.
\end{array}
$$
This implies that the spectral abcissa satisfies,
$$
\varsigma(A-RQ_{\infty})=-\lambda^h_1<0.
$$

Observe that in general situations, even though $(A-RQ_{\infty})$
and $(A-RQ_{\infty})^{\prime}$ have the same eigenvalues we have
$$
(A-RQ_{\infty})^{\prime}=-((A^{\prime})^2+SR)^{1/2}\not=-(A^2+RS)^{1/2}=(A-RQ_{\infty}).
$$
Thus even when $\mu(A)<0$ there are situations where 
$$
\varsigma(A-RQ_{\infty})<0<\mu(A-RQ_{\infty}).
$$
For a more thorough discussion on these situations, we refer the reader to section 4.1 in the article~\cite{tugaut-dm-1}.
Using the formula (\ref{ref-Weyls-2-rev-intro}), as well as the definitions of $h_0$ and $\upsilon$ given in \eqref{def-h} and \eqref{def-upsilon-intro}, respectively, we have
$$
\upsilon(dx)~h_0(x)=~\exp{\left(-\frac{1}{2}~x^{\prime}P_{\infty}^{-1}x\right)}\quad \mbox{\rm and}\quad \upsilon(h_0):=\int \upsilon(dx)~h_0(x)=\sqrt{\mbox{\rm det}(2\pi P_{\infty})}.
$$
This implies that
$$
 \BB_h(\upsilon)=\eta_{\infty}
\quad\mbox{\rm and}\quad
 \eta_{\infty}^h=\BB_h(\eta_{\infty})=\BB_{h^2}(\upsilon).
$$
The limiting covariance matrix of the $h$-process is given more explicitly by the formulae
\begin{equation}\label{dif-form-Ph}
P^h_{\infty}:=
\left(P^{-1}_{\infty}+Q_{\infty}\right)^{-1}=\frac{1}{2}\left(P_{\infty}^{-1} + R^{-1}A\right)^{-1}=\frac{1}{2}R(A^2+RS)^{-1/2},
\end{equation}
where we have used \eqref{ref-Weyls-2-rev-intro} and \eqref{ref-Weyls-2-intro}. Combining this with the second equality in \eqref{ref-Weyls-2-rev-intro}, we obtain
\begin{equation}\label{dif-form-Ph-2}
Q_{\infty}-R^{-1}A=R^{-1}(A^2+RS)^{1/2}=(2P^h_{\infty})^{-1}=\frac{1}{2}\left(P^{-1}_{\infty}+Q_{\infty}\right),
\end{equation}
which implies that
\begin{equation}\label{h02-normalizing}
(\ref{dif-form-Ph-2})\Longrightarrow
\upsilon(h_0^2)={(2\pi)^{r/2}}/{\sqrt{\mbox{\rm det}(P^{-1}_{\infty}+Q_{\infty})}}.
\end{equation}

Finally notice that
\begin{eqnarray}
AR= RA^{\prime}&\Longleftrightarrow &(A-RQ_{\infty})P^h_{\infty}+P^h_{\infty}(A-RQ_{\infty})^{\prime}+R=0\label{sylvester-ref}.
\end{eqnarray}
Thus, our condition ensures the reversibility property (\ref{intro-master}) of the $h$-process.

\begin{rem}
Assume that $S>0$ and $SA=A^{\prime}S$. In this situation, 
the fixed point matrices $(P_{\infty}^-,P_{\infty})$ are given by (cf.~\cite{bp-21})
\begin{equation}
P_{\infty}^-=AS^{-1}-(A^2+RS)^{1/2}S^{-1}< 0<
 P_{\infty}=AS^{-1}+(A^2+RS)^{1/2}S^{-1}. 
 \label{formula-example-Pinfty-0}
 \end{equation}
Thus, whenever $S,R>0$ and $SAS^{-1}=A^{\prime}=R^{-1}AR$ we have
\begin{equation}\label{both-cond}
(A- RQ_{\infty})=-\left(A^2+RS\right)^{1/2}=(A-P_{\infty}S)\quad\mbox{\rm and}\quad
RQ_{\infty}=P_{\infty}S.
\end{equation}
Using (\ref{Delta-P-Q}) we also have
$$
\begin{array}{l}
\displaystyle (P^h_{\infty})^{-1}=
P_{\infty}^{-1}+Q_{\infty}=2\left(P_{\infty}^{-1}+R^{-1}A\right)=2R^{-1}\left(A^2+RS\right)^{1/2}=-2R^{-1}(A-P_{\infty}S)\\
\\
\displaystyle\Longrightarrow\mbox{det}\left(P_{\infty}^{-1}+Q_{\infty}\right)~=2^r~{\sqrt{\vert\mbox{det}\left(A^2+RS\right)\vert}}/{\mbox{det}(R)}.
\end{array}$$
This implies that
$$
\begin{array}{l}
\displaystyle
RQ_{\infty}=P_{\infty}S=RP_{\infty}^{-1}+2A=A+\left(A^2+RS\right)^{1/2}.
\end{array}$$
Whenever $S>0$, up to a change of basis, there is no loss of generality to assume that $S=I$. More precisely the matrices $\overline{P}_t:=S^{1/2}P_tS^{1/2}$ satisfy the same Riccati equation as $P_t$ when we replace $(A,R,S)$ by the matrices 
\begin{equation}\label{ref-overline-A}
	(\overline{A},\overline{R},\overline{S}):=(S^{1/2}AS^{-1/2},S^{1/2}RS^{1/2},I).
\end{equation}
\end{rem}}

\subsection{Chebychev-Hermite polynomials} \label{chb-hermite-sec}
Before moving on to the proofs of the spectral theorems, we start with a brief review on multivariate Chebychev-Hermite polynomials.

The generating function for the family of multivariate Chebychev-Hermite polynomials 
$\HH^r_n(x)$ indexed by $n\in \NN^r$  is defined for any $u,x\in\RR^r$ as the convergent series expansion
$$
\SS^r_u(x):=\sum_{n\in \NN^r}\frac{u^n}{n!}~\HH^r_n(x)=\exp{\left(u^{\prime}x-\frac{1}{2}~u^{\prime}u\right)}
$$
with the multiple index notation
$$
n=(n_1,\ldots,n_r)\in \NN^r\quad u=(u_1,\ldots,u_r)\in\RR^r \Longrightarrow u^n:=u_1^{n_1}\times\ldots\times u_r^{n_r}.
$$
Recall that $\frac{1}{\sqrt{n!}}\HH^r_n$ forms an orthonormal basis of the Hilbert space $\LL_2(\nu)$, where $\nu=\Na(0,I)$ stands for the centered Gaussian measure on $\RR^r$ with unit covariance.
Observe that (\ref{ref-Weyls-2-rev-intro}) implies
\begin{eqnarray}
&&(A-RQ_{\infty})=-(A^2+RS)^{1/2}=-\frac{1}{2}~R(P^h_{\infty})^{-1}\nonumber\\
\displaystyle&\Longrightarrow&\Lambda^h=(P^{h}_{\infty})^{-1/2}(A-RQ_{\infty})(P^{h}_{\infty})^{1/2}=-\frac{1}{2}(P^{h}_{\infty})^{-1/2}R(P^{h}_{\infty})^{-1/2}<0.
\label{commutation-ref}
\end{eqnarray}
In addition, we have
\begin{equation}\label{spec-Lambda-h}
\mbox{\rm Spec}(\Lambda^h)=\left\{-\vert \lambda\vert^{1/2}~:~\lambda\in\mbox{\rm Spec}(A^2+RS) \right\}\subset \RR_-.
\end{equation}
\begin{defi}\label{def-Z}
 We denote by $z_i$ an orthonormal eigenvector of the  matrix $\Lambda^h$ associated with an eigenvalue $\lambda_i(\Lambda^h):=-\lambda_i^h<0$, with $i\in\{1,\ldots,r\}$, and we set 
$$
\Za:=\left(z_1,\ldots,z_r\right)\Longrightarrow \Za^{\prime}\Za=I
\quad \mbox{\rm and}\quad
\Ea^h_t(Q_{\infty}):= \exp{\left(\Lambda^h\,t\right)}.
$$
\end{defi}
\begin{lem}

For any $t\geq 0$ we have
$$
I-(P_{\infty}^h)^{-1/2}P^h_t(P_{\infty}^h)^{-1/2}=
 \Ea^h_t(Q_{\infty})^2\quad \mbox{and}\quad
 \Ea^h_t(Q_{\infty})^{\prime}= \Ea^h_t(Q_{\infty}).
$$
\end{lem}
\proof
{\color{black}Observe that
\textcolor{black}{
\begin{eqnarray*}
AR=RA^{\prime}&\Longrightarrow &R^{-1}(A-RQ_{\infty})=(A-RQ_{\infty})^{\prime}R^{-1}\\
&\Longrightarrow &\forall n\geq 1\quad R^{-1}(A-RQ_{\infty})^nR=\left((A-RQ_{\infty})^{\prime}\right)^n.
\end{eqnarray*}
This yields the formula}
$$
P^h_t=P_{\infty}^h~\left(I-R^{-1}e^{2(A-RQ_{\infty})t}R\right)=P_{\infty}^h~\left(I-e^{2(A-RQ_{\infty})^{\prime}t}\right)
$$}
and therefore
$$
P^h_t=\left(I-e^{2(A-RQ_{\infty})t}\right)P_{\infty}^h.
$$
This implies that
$$
\begin{array}{l}
\displaystyle
(P_{\infty}^h)^{-1/2}P^h_t(P_{\infty}^h)^{-1/2}=I-(P_{\infty}^h)^{-1/2}~e^{2(A-RQ_{\infty})t}(P_{\infty}^h)^{1/2}.
\end{array}
$$
By (\ref{commutation-ref}), we have the commutative property
\begin{equation}\label{h-commutation}
(P^{h}_{\infty})^{-1/2}(A-RQ_{\infty})(P^{h}_{\infty})^{1/2}=(P^{h}_{\infty})^{1/2}(A-RQ_{\infty})^{\prime}(P^{h}_{\infty})^{-1/2},
\end{equation}
which implies that
$$
P_{\infty}^he^{t(A-RQ_{\infty})^{\prime}}=e^{t(A-RQ_{\infty})}P_{\infty}^h\Longleftrightarrow
e^{t(A-RQ_{\infty})}=P_{\infty}^he^{t(A-RQ_{\infty})^{\prime}}(P_{\infty}^h)^{-1}.
$$
Thus, we have
$$
\begin{array}{l}
\displaystyle
(P_{\infty}^h)^{-1/2}~e^{2(A-RQ_{\infty})t}(P_{\infty}^h)^{1/2}\\
\\
\displaystyle=\left(
(P_{\infty}^h)^{-1/2}~e^{(A-RQ_{\infty})t}~(P_{\infty}^h)^{1/2}\right)\left((P_{\infty}^h)^{1/2}e^{t(A-RQ_{\infty})^{\prime}}(P_{\infty}^h)^{-1/2}\right)
\end{array}
$$
and can conclude that
$$
 \Ea^h_t(Q_{\infty})^2= \Ea^h_t(Q_{\infty}) \Ea^h_t(Q_{\infty})^{\prime}.
$$
\cqfd

\begin{lem}\label{lem-S-r-h}
For any $u,x\in\RR^r$ and any $t\geq 0$ we have
\begin{equation}\label{ref-S-r-h}
\EE\left(\SS^r_u\left((P^{h}_{\infty})^{-1/2}X^h_t(x)\right)\right)=\SS^r_{\Ea^h_t(Q_{\infty})u}\left((P^{h}_{\infty})^{-1/2}x\right).
\end{equation}
\end{lem}
\proof
Using the decomposition
$$
(P^{h}_{\infty})^{-1/2}X^{h}_t(x)\stackrel{law}{=}(P^{h}_{\infty})^{-1/2}\widehat{X}^h_t(x)+(P^{h}_{\infty})^{-1/2}(P^{h}_{t})^{1/2} W_1,
$$
where $W_1 \sim \Na(0, I)$, it follows that
$$
\begin{array}{l}
\displaystyle
\log{\EE\left(\exp{\left(u^{\prime}(P^{h}_{\infty})^{-1/2}X^h_t(x)-\frac{1}{2}~u^{\prime}u\right)}\right)}\\
\\
\displaystyle=u^{\prime}(P^{h}_{\infty})^{-1/2}e^{(A-RQ_{\infty})t}x-\frac{1}{2}~u^{\prime}
\left(I-(P^{h}_{\infty})^{-1/2}P^{h}_t(P^{h}_{\infty})^{-1/2}\right)u\\
\\
\displaystyle=u^{\prime}(P^{h}_{\infty})^{-1/2}e^{(A-RQ_{\infty})t}(P^{h}_{\infty})^{1/2}~(P^{h}_{\infty})^{-1/2}x-\frac{1}{2}~u^{\prime}
(P_{\infty}^h)^{-1/2}~e^{2(A-RQ_{\infty})t}(P_{\infty}^h)^{1/2}u.
\end{array}$$
This implies that
\begin{equation}\label{key-expo}
\begin{array}{l}
\displaystyle
\log{\EE\left(\exp{\left(u^{\prime}(P^{h}_{\infty})^{-1/2}X^h_t(x)-\frac{1}{2}~u^{\prime}u\right)}\right)}\\
\\
\displaystyle=\left(\Ea^h_t(Q_{\infty})u\right)^{\prime}~(P^{h}_{\infty})^{-1/2}x-\frac{1}{2}~\left(\Ea^h_t(Q_{\infty})u\right)^{\prime}
\left(\Ea^h_t(Q_{\infty})u\right),
\end{array}\end{equation}
from which the result now follows.
\cqfd

\subsection{Proofs of the spectral theorems}\label{sec-spectral-proofs}
We are now in position to prove the spectral theorems, Theorem~\ref{spectral-theo-intro} and Theorem \ref{theo-Ka-spec}.

\subsubsection{Proof of Theorem~\ref{spectral-theo-intro}}
It suffices to show that
$$
\EE\left(\SS^r_u\left(\Za^\prime({P}^h_{\infty})^{-1/2}X^h_t(x)\right)\right)=\SS^r_{e^{t\overline{\Lambda}^h}u}\left(\Za^{\prime}({P}^h_{\infty})^{-1/2}x\right),
$$
where $\Za$ was defined in Definition \ref{def-Z} and $\overline{\Lambda}^h:=\Za^{\prime}\Lambda^h\Za=\mbox{\rm Diag}\left(-\lambda_1^h,\ldots,-\lambda_r^h\right)$.

To this end, observe that for any $x, u \in \RR^r$, we have
$$
\SS^r_u(\Za^{\prime}x)=\exp{\left(u^{\prime}\Za^{\prime}x-\frac{1}{2}~u^{\prime}\Za^{\prime}\Za u\right)}=\SS^r_{\Za u}(x).
$$
Also note that 
$$
\Za^{\prime}e^{t\Lambda^h}\Za=e^{t\overline{\Lambda}^h} 
\Longrightarrow \Ea_t(Q_\infty)\Za = \Za e^{t\overline{\Lambda}^h}.
$$
Thus, combining these observations with Lemma \ref{lem-S-r-h}, we have
\begin{align*}
\EE\left(\SS^r_u\left(\Za^\prime({P}^h_{\infty})^{-1/2}X^h_t(x)\right)\right) 
&= \EE\left(\SS^r_{\Za u}\left(({P}^h_{\infty})^{-1/2}{X}^h_t(x)\right)\right)\\
&= \SS^r_{\Ea^h_t(Q_{\infty})\Za u}\left(({P}^{h}_{\infty})^{-1/2} x\right)\\
&= \SS^r_{\Za e^{t\overline{\Lambda}^h} u}\left(({P}^{h}_{\infty})^{-1/2} x\right)\\
&= \SS^r_{e^{t\overline{\Lambda}^h} u}\left(\Za^\prime ({P}^{h}_{\infty})^{-1/2} x\right).
\end{align*}
\cqfd

\subsubsection{Proof of Corollary~\ref{Poincare-cor}}
We first prove the estimate \eqref{expo-decays}. Using the decomposition from Theorem \ref{spectral-theo-intro}, it is straightforward to show that for any function $f\in \LL_2(\eta_{\infty}^h)$ we have
\begin{eqnarray*}
\Vert \Ka^h_t(f)-\eta_{\infty}^h(f)\Vert_{2,\eta_{\infty}^h}^2&=&\sum_{n\in\NN^r-\{(0)\}}~e^{-2\lambda_n^ht}~\eta_{\infty}^h(f \varphi^h_n)^2\\
&\leq& e^{-2\lambda_1^ht}
\sum_{n\in\NN^r-\{(0)\}}~~\eta_{\infty}^h(f \varphi^h_n)^2\\
&=& e^{-2\lambda_1^ht}\Vert f -\eta_{\infty}^h(f)\Vert_{2,\eta_{\infty}^h}^2. 
\end{eqnarray*}

\medskip

Now let us prove that this is equivalent to \eqref{poincare-ineq}. For small values of the time parameter $t<(2\lambda^h_1)^{-1}$ an elementary second order Taylor expansion of the exponential function yields
$$
e^{-2\lambda_1^ht}~\Vert f-\eta_{\infty}^h(f)\Vert_{2,\eta_{\infty}^h}^2= 
\left(1-2\lambda_1^ht\right)~\Vert f-\eta_{\infty}^h(f)\Vert_{2,\eta_{\infty}^h}^2+\mbox{\rm o}(t).
$$
In the same vein, we have
$$
\begin{array}{l}
\displaystyle
\Vert \Ka^h_t(f)-\eta_{\infty}^h(f)\Vert_{2,\eta_{\infty}^h}^2\\
\\
\displaystyle=\Vert f-\eta_{\infty}^h(f)\Vert_{2,\eta_{\infty}^h}^2+t~
\eta_{\infty}^h\left(\frac{1}{t}\left(\Ka^h_t(f)-f\right)~\left(\Ka^h_t(f)+f\right)\right)\\
\\
\displaystyle=\Vert f-\eta_{\infty}^h(f)\Vert_{2,\eta_{\infty}^h}^2-2t~E_h(f,f)+\mbox{\rm o}(t).
\end{array}
$$
Using (\ref{expo-decays}) we obtain the Poincar\'e inequality (\ref{poincare-ineq}).

On the other hand, suppose \eqref{poincare-ineq} holds. The Markov transitions $\Ka^h_t$ satisfy
 the Chapman-Kolmogorov evolution equation given in weak form by the formulae
\begin{equation}\label{def-Kolmogorov}
\partial_t\Ka^h_t=\Ka^h_t\La^h=\La^h\Ka^h_t.
\end{equation}
This yields the Dirichlet form equation
$$
\partial_t\, \eta_{\infty}^h(f\,\Ka^h_t(g))_{|t=0}= -E_h(f,g):=\eta_{\infty}^h(f\,\La^h(g)),
$$
and hence
\begin{equation}\label{dirichlet-formulation}
\partial_t\Vert \Ka^h_t(f)-\eta_{\infty}^h(f)\Vert_{2,\eta_{\infty}^h}^2=-2~E_h\left(\Ka^h_t(f),\Ka^h_t(f)\right).
\end{equation}

Combining this with \eqref{poincare-ineq} yields \eqref{expo-decays}. \cqfd

\subsubsection{Proof of Theorem~\ref{theo-Ka-spec}}
Recall that $\Ka_t^h$ and $\Ka_t$ are connected via 
$$
  e^{\lambda_0 t}\Ka_t = \Upsilon_h \circ \Ka_t \circ \Upsilon_h,
$$
where the isometry $\Upsilon_h$ was defined in \eqref{isometry}. Using this and the formula $\eta_\infty^h = \BB_{h_0^2}(\upsilon)$, it is straightforward to show that
$$
  \Ka_t(x, \d y) = \sum_{n\in\NN^r} e^{-\lambda_n t}~\varphi_n(x)~\varphi_n(y)~\upsilon(dy),
$$
where 
$$
\varphi_n(x)=\Upsilon_h \left(\varphi^h_n\right)=\frac{h_0(x)}{\sqrt{\upsilon(h_0^2)}}~\varphi^h_n(x)=\frac{h_0(x)\varphi^h_n(x)}{\sqrt{\upsilon(h_0)\eta_{\infty}(h_0)}}\Longrightarrow
\varphi_0(x)=\frac{h_0(x)}{\sqrt{\upsilon(h_0^2)}}
 $$
 and $\lambda_n = \lambda_0 + \lambda_n^h$. To complete the proof, note that from the definitions of $\eta_\infty$, $h_0$ and $\upsilon$ given in \eqref{cv-eta-infty}, \eqref{def-h} and \eqref{def-upsilon-intro}, respectively, we observe that
 $$
 \eta_{\infty}(h_0)=\left(\mbox{det}\left(I+P_{\infty}Q_{\infty}\right)\right)^{-1/2}\quad \mbox{\rm and}\quad\upsilon(h_0)=\left(\mbox{det}(2\pi P_{\infty})\right)^{1/2}.
 $$
 This yields the formulae
 $$
\upsilon(h_0)\eta_{\infty}(h_0)=
(2\pi)^{r/2}~\left(\mbox{det}\left(P_{\infty}^{-1}+Q_{\infty}\right)\right)^{-1/2},
$$
which ends the proof of the theorem. \cqfd

\subsubsection{Proof of Corollary~\ref{Ka-L2-cor}}
We have the decomposition
\begin{equation}\label{lambda-0-Ka}
e^{\lambda_0t}\Ka_{t}(f)(x)-\frac{h_0(x)}{\eta_{\infty}(h_0)}~ \eta_{\infty}(f)=\sum_{n\in\NN^r-\{0\}} e^{-\lambda^h_nt}~\varphi_n(x)~\upsilon(\varphi_nf).
\end{equation}
The proof of (\ref{expo-decays-Ka}) is now a direct consequence of the formulae $\eta_\infty = \BB_{h_0}(\upsilon)$ and
$$
\sum_{n\in\NN^r-\{0\}} \upsilon(f\varphi_n)^2=\Vert f\Vert_{2,\upsilon}^2-\upsilon(\varphi_0f)^2=\upsilon(f^2)-\left(\upsilon\left(\frac{h_0f}{\sqrt{\upsilon(h_0^2)}}\right)\right)^2.
$$
\cqfd

\section{Appendix}
Here we provide the proofs of the assertions \eqref{de-Bruijn}-\eqref{free-energy-decays} stated in section \ref{rev-models-intro}. Thus we assume that the matrices $(A,R,S)$ satisfy the rank condition (\ref{def-contr-obs}) and that we have $R>0$ and $AR=RA^{\prime}$. We prove the four assertions in the order they are stated.

By the density transport formula (\ref{density-transport}) we have
$$
\eta_0^h(dx)=f_0(x)~\eta^h_{\infty}(dx)\Longrightarrow
\partial_t \mbox{\rm Ent}\left(\eta_t^h~|~\eta_{\infty}^h\right)=\int~(1+\log{\Ka_t^h(f_0)})~\La^h(\Ka_t^h(f_0))~d\eta_{\infty}^h.
$$
Applying the integration by parts formula (\ref{ipb-ref}) to $g=1 + \log{\Ka_t^h(f_0)}$ we find the de Bruijn identity \eqref{de-Bruijn}
$$
\begin{array}{l}
\displaystyle
\partial_t \mbox{\rm Ent}\left(\eta_t^h~|~\eta_{\infty}^h\right)=-\frac{1}{2}~\int~\frac{\Vert \nabla\Ka_t^h(f_0)\Vert^2}{\Ka_t^h(f_0)}~d\eta_{\infty}^h:=-\frac{1}{2}~\Ja\left(\eta_t^h~|~\eta_{\infty}^h\right).
\end{array}$$

\bigskip

Next, observe that using (\ref{Xh-mean-cov}) we obtain
$$
\displaystyle \nabla \widehat{X}^h_t(x)=\exp{\left((A-RQ_{\infty})^{\prime}t\right)}.
$$
Also note that
$$
\nabla\Ka_t^h(f_0)(x)=\nabla\EE(f_0(\widehat{X}^h_t(x)))=
\EE(\nabla \widehat{X}^h_t(x)(\nabla f_0)(\widehat{X}^h_t(x)).
$$
This yields the commutative property
$$
\begin{array}{l}
\displaystyle \nabla\Ka_t^h(f_0)=e^{(A-RQ_{\infty})^{\prime}t}\,\Ka_t^h(\nabla f_0)\\
\\
\displaystyle\Longrightarrow
\Ja\left(\eta_t^h~|~\eta_{\infty}^h\right)=\int \frac{\Ka_t^h(\nabla f_0)(x)}{\sqrt{\Ka_t^h( f_0)(x)}}^{\prime}
 e^{(A-RQ_{\infty})t}e^{(A-RQ_{\infty})^{\prime}t}\frac{\Ka_t^h(\nabla f_0)(x)}{\sqrt{\Ka_t^h(f_0)(x)}}~\eta_{\infty}^h(dx).
 \end{array}
$$
Applying Cauchy Schwartz  inequality we find that
\begin{eqnarray*}
\Ja\left(\eta_t^h~|~\eta_{\infty}^h\right)&\leq &\Vert  e^{(A-RQ_{\infty})t}\Vert^2~\int~\eta_{\infty}^h(dx)~\frac{\Vert \Ka_t^h(\sqrt{f_0}~(\nabla f_0/\sqrt{f_0}))(x)\Vert^2}{\Ka_t^h( f_0)(x)}\\
&\leq &\Vert  e^{(A-RQ_{\infty})t}\Vert^2~\int~\eta_{\infty}^h(dx)~\Ka_t^h(\Vert \nabla f_0\Vert^2/f_0)(x)\\
&=&\Vert  e^{(A-RQ_{\infty})t}\Vert^2~\eta_{\infty}^h(\Vert \nabla f_0\Vert^2/f_0).
\end{eqnarray*}
This yields the Fisher information exponential decays \eqref{exp-decay-Fisher}
$$
\Ja\left(\eta_t^h~|~\eta_{\infty}^h\right)\leq \Vert  \exp{((A-RQ_{\infty})t)}\Vert^2~\Ja\left(\eta_0^h~|~\eta_{\infty}^h\right)\longrightarrow_{t\rightarrow\infty} 0.
$$

\medskip

Applying the de Bruijn identity we have
$$
\begin{array}{l}
\displaystyle
\mbox{\rm Ent}\left(\eta_0^h~|~\eta_{\infty}^h\right)=-
\int_0^{\infty}\partial_s \mbox{\rm Ent}\left(\eta_s^h~|~\eta_{\infty}^h\right)~ds=\frac{1}{2}\int_0^{\infty}~\Ja\left(\eta_s^h~|~\eta_{\infty}^h\right)~ds,
\end{array}$$
which yields the log-Sobolev inequality \eqref{log-Sobolev}
$$
\mbox{\rm Ent}\left(\eta_0^h~|~\eta_{\infty}^h\right)\leq \left(\frac{1}{2}\int_0^{\infty}\Vert  \exp{((A-RQ_{\infty})s)}\Vert^2~ds\right)~\Ja\left(\eta_0^h~|~\eta_{\infty}^h\right).$$

Finally, applying the Log-Sobolev inequality to $\eta^h_t$, the
 de Bruijn identity now yields the free energy exponential decays \eqref{free-energy-decays}
$$
\begin{array}{l}
\displaystyle
\partial_t \mbox{\rm Ent}\left(\eta_t^h~|~\eta_{\infty}^h\right)=-\frac{1}{2}~\Ja\left(\eta_t^h~|~\eta_{\infty}^h\right)
\leq -\left(\int_0^{\infty}\Vert  e^{(A-RQ_{\infty})s}\Vert^2~ds\right)^{-1}~\mbox{\rm Ent}\left(\eta_t^h~|~\eta_{\infty}^h\right).
\end{array}
$$

\end{document}